\documentclass[12pt,a4paper,onecolumn]{article}
\usepackage{textcomp}
\usepackage{tikz}
\usepackage{amssymb}
\usepackage{caption}
\usepackage{color}
\usepackage{fancybox}
\usepackage{verbatim}
\usepackage{fancybox}
\usepackage{titlesec} 
\usepackage{mathtools}
\usepackage{amsmath}
\usepackage{esvect}
\usepackage{cite}
\usepackage{setspace}
\usepackage{diagbox}
\numberwithin{equation}{section}
\usepackage[hidelinks]{hyperref}
\usepackage{mathrsfs}
\usetikzlibrary{positioning,fit,calc}
\usepackage{array}
\usepackage{bm}
\usepackage{multirow}
\usepackage{hhline}
\usepackage[super]{nth}
\usetikzlibrary{positioning,decorations.pathreplacing}
\usepackage{subcaption} 
\usetikzlibrary{decorations.pathmorphing}
\usetikzlibrary{arrows}
\usepackage{fixltx2e}
\usepackage{ragged2e}
\usepackage[super]{nth}

\usepackage{mathrsfs}
\usetikzlibrary{positioning,fit,calc}
\usepackage[a4paper,left=1in,top=1in,right=1in,bottom=1in,nohead]{geometry} 
\usepackage[hidelinks]{hyperref}
\usepackage{array}
\usepackage{bm}
\usepackage[symbol]{footmisc} 
\numberwithin{equation}{section}
\usepackage{amsthm}
\def\correspondingauthor{\footnote{Corresponding author. Email: williewong088@gmail.com.}}

\tikzset{block/.style={draw,thick,text width=2cm,minimum height=1cm,align=center},
         line/.style={-latex}}
\newcolumntype{P}[1]{>{\centering\arraybackslash}m{#1}} 

\titleformat{\section}[block]{\large\scshape\bfseries}{\thesection.}{1em}{} 
\titleformat{\subsection}[block]{\bfseries}{\thesubsection.}{1em}{} 

\newtheorem{thm}{Theorem}[section]
\newtheorem{ppn}[thm]{Proposition}
\newtheorem{cor}[thm]{Corollary}
\newtheorem{lem}[thm]{Lemma}

\theoremstyle{definition}
\newtheorem{defn}[thm]{Definition}

\newtheorem{eg}[thm]{Example}

\newtheorem{prob}[thm]{Problem}

\begin{document}
\pagenumbering{arabic}
\begin{center}
    \textbf{\Large Arc reversals of cycles in\\orientations of $\bm{G}$ vertex-multiplications}
\vspace{0.1 in} 
    \\{\large H.W. Willie Wong\correspondingauthor{}, E.G. Tay}
\vspace{0.1 in} 
\\National Institute of Education\\Nanyang Technological University\\Singapore
\end{center}

\begin{abstract}
\noindent Ryser proved that any two tournaments with the same score sequence are $C_3$-equivalent while Beineke and Moon proved the $C_4$-equivalence for any two bipartite tournaments with the same score lists. In this paper, we extend these results to orientations of $G$ vertex-multiplications. We focus on two main areas, namely orientations with the same score list and with score-list parity. Our main tools are extensions of the refinement technique, directed difference graph and a reduction lemma.
\end{abstract}
\section{Introduction}
In this paper, let $G$ be a connected graph with vertex set $V(G)=\{1,2,\ldots, n\}$ and edge set $E(G)$. We consider only graphs with no loops or parallel edges. An edge $e\in E(G)$ is a \textit{bridge} if $G-e$ is disconnected. For any $v\in V(G)$, the \textit{degree} of $v$, denoted by $deg_G(v)$, is the number of edges incident to it. For any vertices $v,x\in V(G)$, the $\textit{distance}$ from $v$ to $x$, $d_G(v,x)$, is defined as the length of a shortest path (i.e., number of edges) from $v$ to $x$. The subgraph of $G$ induced by the set of vertices $V\subseteq V(G)$ (set of edges $E\subseteq E(G)$ resp.) is denoted by $\langle V\rangle_G$ ($\langle E\rangle_G$ resp.). The definitions of distance and induced subdigraph are analogous for digraphs; and we refer the reader to \cite{BJJ GG} and \cite{JWM3} for any undefined terminology. We denote the complete graph of order $n$ by $K_n$ and the complete bipartite graph with partite sets $X$ and $Y$ by $K(X,Y)$. A \textit{chord} of a cycle $C$ is an edge not in $C$ whose endpoints lie in $C$ and a \textit{chordless} cycle is a cycle that has no chord. A graph is \textit{chordal} if it has no chordless cycle of length at least $4$. For any graph $G$, let $\mathcal{L}^*_G$ ($\mathcal{L}_G$ resp.) denote the set of lengths of all cycles (chordless cycles resp.) in $G$, i.e.,
\begin{align*}
\mathcal{L}^*_G=\{l\mid G \text{ contains an $l$-cycle}\} \text{ and }\mathcal{L}_G=\{l\mid G \text{ contains a chordless $l$-cycle}\}.
\end{align*}
For instance, $\mathcal{L}^*_{K_n}=\{3,4,\ldots, n\}$ and $\mathcal{L}_{K_n}=\{3\}$.
\noindent\par Let $D$ be a digraph with vertex set $V(D)$ and arc set $A(D)$. If $uv\in A(D)$, where $u,v\in V(D)$, then we say $u$ \textit{dominates} $v$ and denote this by $u\rightarrow v$. The \textit{outset} and \textit{inset} of a vertex $v\in V(D)$ are defined to be $O_D(v)=\{x\in V(D)\mid v\rightarrow x\}$ and $I_D(v)=\{y\in V(D)\mid y\rightarrow v\}$ respectively. The \textit{score} $s_v$ or \textit{outdegree} $od_D(v)$ of a vertex $v\in V(D)$ is defined by $s_v=od_D(v)=|O_D(v)|$. That is, we shall freely interchange between the two notations, $s_v$ and $od_D(v)$. The \textit{co-score} or \textit{indegree} $id_D(v)$ of a vertex $v\in V(D)$ is defined by $id_D(v)=|I_D(v)|$. The \textit{underlying graph} of $D$, denoted by $\mathrm{U}(D)$, is the undirected graph with the same vertex set and all arcs replaced by undirected edges. If there is no ambiguity, we shall omit the subscript for the above notation.
\indent\par An $\textit{orientation}$ $D$ of a graph $G$ is a digraph obtained from $G$ by assigning a direction to every edge in $E(G)$. Therefore, a tournament (multipartite or $n$-partite tournament resp.) is an orientation of the complete graph (complete $n$-partite graph resp.). For $k\ge 3$, we denote the $k$-dicycle (i.e., directed cycle of length $k$) by $C_k$ and the collection of all orientations of a $k$-cycle by $\mathscr{C}_k$. The \textit{converse} $\tilde{D}$ of any digraph $D$ is the digraph with $V(\tilde{D})=V(D)$ and $A(\tilde{D})=\{uv\mid vu\in A(D)$\}. We say that $D$ is \textit{self-converse} if $D$ is isomorphic to $\tilde{D}$, i.e., $D\cong \tilde{D}$. Also, we say a collection $\mathscr{F}$ of digraphs is \textit{self-converse} if $\mathscr{F}=\{\tilde{F}\mid F\in \mathscr{F}\}$. For instance, it is easy to verify that $D$ is self-converse for all $D\in \mathscr{C}_3\cup \mathscr{C}_4$ and $\mathscr{C}_k$ is a self-converse collection for all $k\ge 3$. We will make use of the following powerful principle quoted from Harary \cite{FH} and use phrases such as ``by duality" or ``the dual of" to refer to this principle.
\begin{justify}
``\textbf{Principle of Directional Duality} For each theorem about digraphs, there is a corresponding theorem obtained by replacing every concept by its converse."
\end{justify} 
\indent\par One of the early problems concerning arc reversals is: Given two tournaments of the same order, is it possible to obtain one from the other by a sequence of a prescribed type of arc reversals? In 1964, Ryser \cite{HJR} gave an affirmative answer of using $C_3$ for any two tournaments with the same score sequence. Later, Waldrop \cite{CW1} gave an independent proof and further established two results in which ``$C_3$" is replaced by ``$C_4$" and by ``$C_5$".  The investigation was extended to bipartite tournaments by Beineke and Moon \cite{LWB JWM}.
\begin{thm} (Ryser \cite{HJR}, Waldrop \cite{CW1}) \label{thmA9.1.1}
Two tournaments have the same score sequence if and only if they are $C_3$-equivalent.
\end{thm}
\begin{thm}(Beineke and Moon \cite{LWB JWM}) \label{thmA9.1.2}
Two bipartite tournaments have the same score lists if and only if they are $C_4$-equivalent.
\end{thm}

\noindent\par On a related note, Reid \cite{KBR} also gave an affirmative answer via $k$-dipath reversals, for all $1\le k \le n-1$, for any two tournaments of order $n$. Some generalisations of these results were obtained by Thomassen \cite{CT}.
\indent\par In this paper, we extend Ryser's and Beineke and Moon's theorems to a large family of graphs known as the $G$ vertex-multiplications. Our work on $G$ vertex-multiplications is still in the early stages and we hope this work illustrates the potential of $G$ vertex-multiplications as new research grounds for arc reversals. We focus on two areas, namely orientations of $G$ vertex-multiplications with the same score list and with score-list parity in Sections 3 and 4 respectively. Our main tools are extensions of the refinement technique and directed difference graph, following Waldrop \cite{CW1} and Thomassen \cite{CT}, and a reduction lemma; we derive them in Section 2. Many of these results have their analogues in tournaments and bipartite tournaments, and references will be provided as we proceed. It will also be notable that the results in bipartite tournaments generalise particularly well to orientations of tree vertex-multiplications.
\noindent\par In 2000, Koh and Tay \cite{KKM TEG 8} introduced the $G$ vertex-multiplications and extended the results on the orientation number of complete $n$-partite graphs. This direction of research was pursued by several others with emphasis on vertex-multiplications of trees and of cycles (see Koh and Tay \cite{KKM TEG 11}, Ng and Koh \cite{NKL KKM} and Wong and Tay \cite{WHW TEG 3A, WHW TEG 6A, WHW TEG 7A, WHW TEG 8A}). 
\noindent\par Let $G$ be a given connected graph. For any sequence of $n$ positive integers $(p_i)$, a $G$ vertex-multiplication, denoted by $G(p_1, p_2,\ldots, p_n)$, is the graph with vertex set $V^*=\bigcup_{i=1}^n{V_i}$ and edge set $E^*$, where $V_i$'s are pairwise disjoint sets with $|V_i|=p_i$, for $i=1,2,\ldots,n$; and for any $u,v\in V^*$, $uv\in E^*$ if and only if $u\in V_i$ and $v\in V_j$ for some $i,j\in \{1,2,\ldots, n\}$ with $i\neq j$ such that $ij\in E(G)$. For instance, if $G\cong K_n$, then the graph $G(p_1, p_2,\ldots, p_n)$ is a complete $n$-partite graph with partite sizes $p_1, p_2,\ldots, p_n$. For convenience, we denote $|V^*|=\sum\limits_{i=1}^{n}{p_i}=N$. Also, we say $G$ is a \textit{parent graph} of $G(p_1, p_2,\ldots, p_n)$.
\indent\par Let $D$ be an orientation of $G(p_1, p_2, \ldots, p_n)$ with the score list $S=(s_1, s_2, \ldots, s_N)$. We shall always assume that $S$ is a list of non-negative integers in non-decreasing order, i.e., $0\le s_1\le s_2\le \ldots\le s_N$, and the vertices are labelled such that $s_i$ corresponds to the vertex $v_i$ for $i=1,2,\ldots, N$. The same holds for any other orientation $D'$ of $G(p_1, p_2, \ldots, p_n)$ with score list $S'=(s'_i)$ and a customised labelling $v'_i$ of vertices. We will loosely use the two denotations of a vertex to refer to its score, for example, if $v_i=w$, then $s_i=s_w$. Let $F$ be a subdigraph of $D$. The \textit{underlying partite} of $F$ is the set $\mathrm{UP}_G(F)=\{i\in V(G) \mid u\in V_i \text{ for some }u\in V(F)\}$. Let $E\subseteq E(G)$ such that $ij\in E$ if and only if $uv\in A(F)$ for some $u\in V_i, v\in V_j$ or $u\in V_j, v\in V_i$. The \textit{underlying partite graph} of $F$ is the graph $\mathrm{UPG}_G(F)=\langle E\rangle_G.$ So, $\mathrm{UPG}_G(F)$ is a spanning subgraph of $\langle \mathrm{UP}_G(F)\rangle_G$. The $G$ vertex-multiplications is a natural generalisation of complete multipartite graphs. Analogous to multipartite tournaments, one can interpret $D$ as the result of competition among the individuals on $n$ teams or comparisons between the items in $n$ sets according to the adjacency of $G$.
\indent\par Lastly, we introduce formally some terminology on arc reversals. For a given $G$, let $\mathscr{D}$ be an arbitrary family of orientations of a $G(p_1, p_2, \ldots, p_n)$ and $\mathscr{F}$ be a collection of asymmetric digraphs. For any $D\in \mathscr{D}$, an \textit{$\mathscr{F}$-reversal} in $D$ is the operation $D\mapsto D_1$ of reversing all the arcs of a copy $F_0$ in $D$ of a digraph $F\in \mathscr{F}$, i.e., $F_0\cong F$, to yield a new orientation $D_1\in\mathscr{D}$. That is, $V(D_1)=V(D)$ and $A(D_1)=(A(D)-A(F_0))\cup A(\tilde{F_0})$. Indeed, $D$ must contain some copy of $F\in\mathscr{F}$ in order for any $\mathscr{F}$-reversal to be performed so we consider only collections of asymmetric digraphs $\mathscr{F}$ for orientations $D$. We denote $D \sim D'$ via \textit{$\mathscr{F}$-reversals} if there is a (possibly null) sequence $D=D_0\mapsto D_1 \mapsto \ldots \mapsto D_k \cong D', k\ge 0$, of $\mathscr{F}$-reversals beginning with $D$ and ending with a digraph isomorphic to $D'$. If $\mathscr{F}$ is a singleton, say $\mathscr{F}=\{F\}$, then we write $\mathscr{F}$-reversal simply as $F$-reversal. It is easy to see that the relation $\sim $ is reflexive and transitive but may not be symmetric. Henceforth, we shall concern ourselves with only self-converse collections $\mathscr{F}$ so that the relation $\sim$ is an equivalence relation on $\mathscr{D}$; we refer to these equivalence classes as $\mathscr{F}$-classes. Furthermore, we say that $D$ and $D'$ are \textit{equivalent via $\mathscr{F}$-reversals} on $\mathscr{D}$, or simply \textit{$\mathscr{F}$-equivalent}, if $D\sim D'$. 

\section{Main tools}
In this section, we derive our main tools for subsequent sections. The first tool is an extension of the refinement technique by Waldrop \cite{CW1} and shares a similar proof. Let $\mathscr{F}_1$ and $ \mathscr{F}_2$ be two self-converse collections of asymmetric digraphs and $\mathscr{D}$ be an arbitrary family of orientations of $G(p_1, p_2, \ldots, p_n)$. We say that $\mathscr{D}$ is \textit{closed under $\mathscr{F}_2$-reversals}, or \textit{simply $\mathscr{F}_2$-closed}, if for any $D\in \mathscr{D}$, an $\mathscr{F}_2$-reversal performed on $D$ yields $D'\in \mathscr{D}$. We say that $\mathscr{F}_1$ \textit{refines} $\mathscr{F}_2$ in the family $\mathscr{D}$ if given any $D\in\mathscr{D}$ and any copy $F_0$ of any digraph $F\in \mathscr{F}_2$, there exists a (possibly null) sequence $D=D_0\mapsto D_1 \mapsto \ldots \mapsto D_k=D', k\ge 0$, of $\mathscr{F}_1$-reversals which \textit{exactly reverses} $F_0$. That is, $V(D')=V(D)$ and $A(D')=(A(D)-A(F_0))\cup A(\tilde{F_0})$. Note that the definition requires the equality $D_k=D'$ and not merely isomorphism. For example, it is easy to prove that $C_3$ refines $C_k$, $k\ge 3$, in the family of all tournaments (see Corollary 1.4.6 in \cite{CW1} for details).

\begin{lem} \textbf{(Refinement)}
Let $G$ be a graph and suppose $\mathscr{F}_1, \mathscr{F}_2$ are self-converse collections of asymmetric digraphs such that $\mathscr{F}_1$ refines $\mathscr{F}_2$ in an arbitrary family $\mathscr{D}$ of orientations of $G(p_1, p_2, \ldots, p_n)$. If $\mathscr{D}$ is $\mathscr{F}_2$-closed and $D, D'\in \mathscr{D}$ are $\mathscr{F}_2$-equivalent, then $D$ and $D'$ are $\mathscr{F}_1$-equivalent. Hence, every $\mathscr{F}_2$-class is contained in an $\mathscr{F}_1$-class.
\end{lem}
\noindent\textit{Proof}: Since $D, D'\in \mathscr{D}$ are $\mathscr{F}_2$-equivalent, there exists a (possibly null) sequence $D=D_0\mapsto D_1 \mapsto \ldots \mapsto D_k\cong D', k\ge 0$, of $\mathscr{F}_2$-reversals. And, $\mathscr{F}_1$ refines $\mathscr{F}_2$ implies for each $0\le i\le k-1$, there exists a sequence $D_i=D^0_i\mapsto D^1_i \mapsto\ldots \mapsto D^{j_i}_i=D_{i+1}$ of $\mathscr{F}_1$-reversals. Adjoining all such sequences in order forms a sequence of $\mathscr{F}_1$-reversals that obtains $D'$ from $D$.
\qed

\indent\par The directed difference graph is another useful tool by Waldrop \cite{CW1} and Thomassen \cite{CT}. We prove an extension prescribed for $G$ vertex-multiplications. Recall that an \textit{automorphism} of $G$ is a bijection $f:V(G)\mapsto V(G)$ such that $uv\in E(G)$ if and only if $f(u)f(v)\in E(G)$.
\begin{defn}\label{defnA9.2.2}
Let $G$ be a graph and $D$ and $D'$ be two orientations of $G(p_1, p_2, \ldots, p_n)$. If $f:V(\mathrm{U}(D))\mapsto V(\mathrm{U}(D'))$ is a automorphism, then the \textit{directed difference graph} $\mathcal{D}=DDG(f;D, D')$ of $D$ and $D'$ is defined as $V(\mathcal{D})=V(D)$ and for any $u,v\in V(\mathcal{D})$, 
\begin{align}
u\rightarrow v \text{ in } \mathcal{D}\iff u\rightarrow v \text{ in } D \text{ and } f(v)\rightarrow f(u) \text{ in } D'. \label{eqA9.2.1}
\end{align}
\end{defn}

\begin{lem}\label{lemA9.2.3}
Let $G$ be a graph and $f$, $D, D'$ and $\mathcal{D}$ be given as in Definition \ref{defnA9.2.2}. Then for any $v\in V(\mathcal{D})$, $od_\mathcal{D}(v)-id_\mathcal{D}(v)=od_{D}(v)-od_{D'}(f(v))=id_{D'}(f(v))-id_{D}(v)$.
\end{lem}
\noindent\textit{Proof}: We prove the first equality and the second equality follows from the fact $od_{D}(v)+id_{D}(v)=od_{D'}(f(v))+id_{D'}(f(v))$. For any $v\in V(\mathcal{D})$, note by (\ref{eqA9.2.1}) that
\begin{align*}
O_{D}(v)&=\{ u\in O_{D}(v)\mid f(u)\in O_{D'}(f(v))\}\sqcup \{u\in O_{D}(v)\mid f(u)\in I_{D'}(f(v))\}\\
&=\{ u\in O_{D}(v)\mid f(u)\in O_{D'}(f(v))\}\sqcup O_\mathcal{D}(v)\text{ and}\\
O_{D'}(f(v))&=\{w\in O_{D'}(f(v)) \mid f^{-1}(w)\in O_{D}(v)\}\sqcup \{w\in O_{D'}(f(v)) \mid f^{-1}(w) \in I_{D}(v)\}\\
&=\{w\in O_{D'}(f(v)) \mid f^{-1}(w)\in O_{D}(v)\}\sqcup I_\mathcal{D}(v),
\end{align*}
where $\sqcup$ denotes disjoint union. Subtracting the cardinalities of the latter from the former gives the first equality.
\qed

\indent\par The next lemma is an extension of a technique sometimes known as cycle stacking.
\begin{lem} \textbf{(Reduction)}
Let $G$ be a graph and $D$ be an orientation of $G(p_1, p_2, \ldots, p_n)$. Suppose $Z=u_1 u_2 \ldots u_k u_1$ is a $k$-dicycle, $k\ge 5$, in $D$ such that for some $1\le i<j\le k$, $u_i$ and $u_j$ are in the same partite set. Then, a $Z$-reversal is equivalent to reversing the arcs of the dicycles 
\begin{align}
Z'=u_{i-1} u_{j} u_{j+1} \ldots u_k u_1 \ldots u_{i-1}\text{ and }Z''=u_{i-1} u_i \ldots u_j u_{i-1}\label{eqA9.2.2}
\end{align}
in some order. Furthermore$^*$, if $H=\mathrm{UPG}(u_i u_{i+1} \ldots u_j)$ is a tree in $G$, then a $Z''$-reversal is equivalent to some sequence of $C_4$-reversals.
\end{lem}
\noindent\textit{Proof}: Consider the cases $u_{i-1}\rightarrow u_j$ or $u_j\rightarrow u_{i-1}$. If the former holds, then a $Z$-reversal is equivalent to a $Z'$-reversal followed by a $Z''$-reversal. If the latter holds, then swapping the order of $Z'$-reversal and $Z''$-reversal achieves the same desired result.
\indent\par For the ``Furthermore" part, note that since $H$ is bipartite, it follows that $j\equiv i \pmod{2}$. We proceed by induction on the difference $j-i$. In the base case of $j=i+2$, $Z''\cong C_4$. Next, assume $j=i+2r$ for some $r\ge 2$ and the partite set $V_1$ contains $u_i$ and $u_j$. Since $H$ is a tree, there exists some vertex $w\in V(H)$ with the largest distance $d_H(1,w)$ and $u_{x+1}\in V_w\cap V(Z)$ for some $i<x+1< j$. Then, $u_x$ and $u_{x+2}$ are in the same partite set. By considerations similar to the base case with $x$ in place of $i$, we deduce that a $Z''$-reversal is equivalent to reversing the arcs of the dicycle $Z_x=u_{i-1} u_i\ldots u_{x-1} u_{x+2} u_{x+3}\ldots u_j u_{i-1}$ and a $4$-dicycle $u_{x-1} u_x u_{x+1} u_{x+2} u_{x-1}$ in some order. By induction hypothesis, a $Z_x$-reversal, and thus a $Z''$-reversal, is feasible via $C_4$-reversals.
\qed

\indent\par We omit the proof of the next lemma in view of its similarity with the previous lemma.
\begin{lem} \textbf{(Generalised Reduction)}
Let $G$ be a graph and $D$ be an orientation of $G(p_1, p_2, \ldots, p_n)$. Suppose $F$ is an orientation of the $k$-cycle $u_1 u_2 \ldots u_k u_1$, $k\ge 5$, in $D$ such that for some $1\le i<j\le k$, $u_i$ and $u_j$ are in the same partite set. Also, denote the arc joining $u_{i-1}$ and $u_j$ as $a_{i-1,j}$. Then, an $F$-reversal is equivalent to reversing the arcs of the orientations 
\begin{align*}
F'=\langle \{u_{j}, u_{j+1}, \ldots, u_k, u_1, \ldots, u_{i-1}\}\rangle_F\cup \{a_{i-1,j}\} \text{ and }
F''=\langle \{u_{i-1}, u_i, \ldots, u_j\}\rangle_F\cup \{\tilde{a}_{i-1,j}\}
\end{align*}
of cycles. Furthermore$^*$, if $F^*=\langle \{u_i, u_{i+1},\ldots, u_j\}\rangle_F$ and $H=\mathrm{UPG}(F^*)$ is a tree in $G$, then an $F''$-reversal is equivalent to some sequence of $\mathscr{C}_4$-reversals.
\end{lem}
{\color{blue} Proof omitted from journal paper
\\\noindent\textit{Proof}: Clearly, an $F$-reversal is equivalent to an $F'$-reversal followed by an $F''$-reversal. For the ``Furthermore" part, note that since $H$ is bipartite, it follows that $j\equiv i \pmod{2}$. We proceed by induction on the difference $j-i$. In the base case of $j=i+2$, $F''\in \mathscr{C}_4$. Next, assume $j=i+2r$ for some $r\ge 2$ and the partite set $V_1$ contains $u_i$ and $u_j$. Since $H$ is a tree, there exists some vertex $w\in V(H)$ with the largest distance $d_H(1,w)$ and $u_{x+1}\in V_w\cap V(F)$ for some $i<x+1< j$. Then, $u_x$ and $u_{x+2}$ are in the same partite set and denote the arc joining $u_{x-1}$ and $u_{x+2}$ in $D$ as $a_{x-1,x+2}$. By considerations similar to the base case with $x$ in place of $i$, we deduce that an $F''$-reversal is equivalent to reversing the arcs of the orientations $F_x=\langle \{u_{i-1}, u_i, \ldots, u_{x-1}, u_{x+2}, u_{x+3}, \ldots, u_j\}\rangle_{F''} \cup\{a_{x-1,x+2}\}$ and $\langle \{u_{x-1}, u_x, u_{x+1}, u_{x+2}\}\rangle_{F''}\cup\{\tilde{a}_{x-1,x+2}\}$ of cycles; note that the latter belongs to $\mathscr{C}_4$. By induction hypothesis, an $F_x$-reversal, and thus an $F''$-reversal, is feasible via $\mathscr{C}_4$-reversals.
\qed
}

~\noindent\par For brevity, we shall refer to the ``Furthermore" part of the Reduction Lemma (Generalised Reduction Lemma resp.) as Reduction$^*$ Lemma (Generalised Reduction* Lemma resp.).

\section{Orientations with the same score list}
Motivated by Theorems \ref{thmA9.1.1} and \ref{thmA9.1.2}, we study orientations of $G$ vertex-multiplications with the same score list in this section. Let $D$ and $D'$ be two orientations of $G(p_1, p_2, \ldots, p_n)$ with score lists $S$ and $S'$ respectively. We say $D$ and $D'$ have the \textit{same score list} if there exists some automorphism $f$ of $G(p_1, p_2, \ldots, p_n)$ such that $od_{D}(v)=od_{D'}(f(v))$ for all $v\in V(D)$. We should mention that in Theorem \ref{thmA9.1.2} the ``score lists" of a bipartite tournament $D$ ($D'$ resp.) refer to two sequences $S_1, S_2$ ($S'_1, S'_2$ resp.) of scores in non-decreasing order and corresponding to its two partite sets; furthermore, $D$ and $D'$ are said to have the same score list if $\{S_1, S_2\}=\{S'_1, S'_2\}$. These two definitions of ``same score list(s)" are equivalent in the context of $G\cong K_2$.
\indent\par Our first two results relax $C_4$-reversals to arc reversals of any dicycles (even dicycles resp.) to extend Theorem \ref{thmA9.1.2} to orientations of $G(p_1, p_2, \ldots, p_n)$ for any graph $G$ (bipartite graph $G$ resp.). Since a $C_k$-reversal does not change the score of any vertex, the sufficiencies of Theorems \ref{thmA9.1.1} and \ref{thmA9.1.2} are trivially true. For the same reason, we omit the sufficiency proofs of the following results.

\begin{thm}\label{thmA9.3.1}
Let $G$ be a graph. Two orientations of $G(p_1,p_2, \ldots, p_n)$ have the same score list if and only if each orientation can be obtained from the other by successively reversing the arcs of $k$-dicycles, $k\ge 3$, whose underlying cycles are pairwise edge-disjoint.
\end{thm}
\noindent\textit{Proof}: $(\Rightarrow)$ Let $D$ and $D'$ be two orientations of $G(p_1, p_2, \ldots, p_n)$ with the same score list. By definition, there exists some automorphism $f$ of $G(p_1, p_2, \ldots, p_n)$ such that $od_{D}(v)=od_{D'}(f(v))$ for all $v\in V(D)$. By Lemma \ref{lemA9.2.3} and denoting $\mathcal{D}=DDG(f;D,D')$,
\begin{align*}
od_\mathcal{D}(v)-id_\mathcal{D}(v)=od_{D}(v)-od_{D'}(f(v))=0
\end{align*} for all $v\in V(\mathcal{D})$. That is, $\mathcal{D}$ is an Eulerian digraph and can be partitioned into dicycles whose underlying cycles are edge-disjoint. Therefore, the result follows.
\qed

\begin{cor}\label{corA9.3.2}
Let $G$ be a bipartite graph. Two orientations of $G(p_1,p_2, \ldots, p_n)$ have the same score list if and only if each orientation can be obtained from the other by successively reversing the arcs of $2k$-dicycles, $k\ge 2$, where the underlying cycles are pairwise edge-disjoint.
\end{cor}
\noindent\textit{Proof}: $(\Rightarrow)$ Invoking Theorem \ref{thmA9.3.1}, it suffices to show that all required $C_k$-reversals are even dicycles, which follows from the fact that $G(p_1, p_2, \ldots, p_n)$ is bipartite.
\qed

\indent\par If the edge-disjoint condition in Theorem \ref{thmA9.3.1} is dropped, one should expect an improvement which require $C_k$-reversals for some (but not all) $k\ge 3$. For the rest of this section, we shall undertake this task, starting with the family of orientations of tree vertex-multiplications.

\begin{ppn}\label{ppnA9.3.3}
Let $T$ be a tree. 
\\(i) In the family of orientations of $T(p_1, p_2, \ldots, p_n)$, $C_4$ refines $\{C_{2k}\mid k\ge 2\}$.
\\(ii) Two orientations of $T(p_1,p_2, \ldots, p_n)$ have the same score list if and only if they are $C_4$-equivalent.
\end{ppn}
\noindent\textit{Proof}: (i) Let $D$ be an orientation of $T(p_1,p_2, \ldots, p_n)$. Since $T(p_1, p_2, \ldots, p_n)$ is bipartite, $D$ contains no odd dicycles. Furthermore, $T$ is a tree implies that any even dicycle $u_1 u_2 \ldots u_{2k} u_1$, $k\ge 2$, in $D$ contains some pair of vertices $u_i$ and $u_j$ from the same partite set. If $k=2$, then we are done. Consider $k\ge 3$. By relabelling the vertices if necessary, we may assume $1\le i<j\le 2k$ and invoke the Reduction$^*$ Lemma (repeatedly) to deduce (i).
\indent\par (ii) $(\Rightarrow)$ By Theorem \ref{thmA9.3.1}, two orientations of $G(p_1,p_2, \ldots, p_n)$ with the same score list are $\{C_k\mid k\ge 3\}$-equivalent. So, they are $C_4$-equivalent by (i), the absence of odd cycles in $T(p_1, p_2, \ldots, p_n)$ and the Refinement Lemma.
\qed

\indent\par It is easy to see that Proposition \ref{ppnA9.3.3}(ii) reduces to Theorem \ref{thmA9.1.2} by setting $G=K_2$. Next, we generalise Proposition \ref{ppnA9.3.3} to any graph $G$.

\begin{thm}\label{thmA9.3.4}
Let $G$ be a graph.
\\(i) In the family of orientations of $G(p_1, p_2, \ldots, p_n)$, $\{C_i\mid i\in \{4\}\cup \mathcal{L}_G\}$ refines $\{C_k \mid k\ge 3\}$.
\\(ii) Two orientations of $G(p_1,p_2, \ldots, p_n)$ have the same score list if and only if they are $\{C_i\mid i\in \{4\}\cup \mathcal{L}_G\}$-equivalent.
\end{thm}
\noindent\textit{Proof}: (i) Let $D$ be an orientation of $G(p_1, p_2, \ldots, p_n)$ containing a $k$-dicycle $Z=u_1 u_2\ldots u_k u_1$, $k\ge 3$. We start by proving a claim.
\\
\\Claim: $\{C_i\mid i\in \{4\}\cup \mathcal{L}^*_G\}$ refines $\{C_k \mid k\ge 3\}$.
\indent\par If $k=3$, then $\mathrm{UPG}(Z)$ is a $3$-cycle in $G$, i.e., $k\in \mathcal{L}^*_G$. It suffices to consider $k\ge 5$. We proceed by induction on the number of pairs $p$ of vertices $u_i, u_j$ of $Z$ in a common partite set. If $p=0$, then $\mathrm{UPG}(Z)$ is a cycle in $G$ and $k\in\mathcal{L}^*_G$. Suppose $p\ge 1$. By relabelling vertices if necessary, we may assume $u_i, u_j\in V(Z)\cap V_w$ for some $1\le i<j\le k$ and $w\in V(G)$; we further assume the difference $j-i$ to be the smallest among all such pairs. If $H=\mathrm{UPG}(u_i u_{i+1}\ldots u_j)$ is a tree in $G$, then by the Reduction$^*$ Lemma, a $Z$-reversal is equivalent to some sequence of $C_4$-reversals and a $Z'$-reversal (as in (\ref{eqA9.2.2})). Since $Z'$ is a dicycle with at most $p-1$ pairs of vertices in a common partite set, the claim follows by induction hypothesis. \indent\par Hence, consider the case where $H$ contains some cycle in $G$. By the Reduction Lemma with $Z'$ and $Z''$ as in (\ref{eqA9.2.2}), a $Z$-reversal is equivalent to a $Z'$-reversal and a $Z''$-reversal in some order. By invoking the Reduction$^*$ Lemma on $Z''$ again, we deduce a $Z''$-reversal is equivalent to reversing the arcs of the dicycle $Z'''=u_{j-1} u_{i} u_{i+1} \ldots u_{j-1}$ and the $4$-dicycle $u_{i-1} u_i u_{j-1} u_j u_{i-1}$ in some order. By minimality of $j-i$, $\mathrm{UPG}(Z''')$ is a cycle in $G$; and since $Z'$ has at most $p-1$ pairs of vertices in a common partite set, the claim follows by induction hypothesis. This completes the proof of the claim.
\\
\indent\par By the above proof, we may assume that $\mathrm{UPG}(Z)$ is a cycle in $G$ and no two vertices of $Z$ belong to the same partite set. It remains to show that a $Z$-reversal is equivalent to some sequence of $\{C_i\mid i\in \{4\}\cup \mathcal{L}_G\}$-reversals. We proceed by induction on the number of chords in $\langle \mathrm{UP}(Z)\rangle_G$. If $k= 4$ or $\langle \mathrm{UP}(Z)\rangle_G$ is chordless (which includes the case $k=3$), then $k\in \{4\}\cup\mathcal{L}_G$ and we are done. So, consider $k\ge 5$ and $\langle \mathrm{UP}(Z)\rangle_G$ has $t$ chords, $t\ge 1$. Let $a_{1,i}$ denote the arc joining $u_1$ and $u_i$ in $D$. WLOG, assume that there exists some $3\le i\le k-1$ such that $\mathrm{UPG}(a_{1,i})$ is a chord of $\langle \mathrm{UP}(Z)\rangle_G$ and $\langle \mathrm{UP}(u_i u_{i+1} \ldots u_k u_1)\rangle_G$ is a chordless cycle in $G$. If $u_i\rightarrow u_1$ in $D$, then a $Z$-reversal is equivalent to reversing the arcs of dicycles $Z_1=u_1 u_2 \ldots u_i u_1$ and $Z_2=u_i u_{i+1} \ldots u_k u_1 u_i$. If $u_1\rightarrow u_i$, then swapping the order of $Z_1$-reversal and $Z_2$-reversal achieves the same desired result. Note that no two vertices of $Z_2$ belong to the same partite set and $\langle \mathrm{UP}(Z_2)\rangle_G$ is a chordless cycle in $G$, i.e., $Z_2\in \{C_i\mid i\in \mathcal{L}_G\}$. Since $\langle \mathrm{UP}(Z_1)\rangle_G$ has at most $t-1$ chords, the result follows by induction hypothesis.
\indent\par (ii) $(\Rightarrow)$ This follows from Theorem \ref{thmA9.3.1}, (i) and the Refinement Lemma.
\qed

\begin{cor}\label{corA9.3.5}
Let $G$ be a chordal graph.
\\(i) In the family of orientations of $G(p_1, p_2, \ldots, p_n)$, $\{C_3, C_4\}$ refines $\{C_k \mid k\ge 3\}$.
\\(ii) Two orientations of $G(p_1,p_2, \ldots, p_n)$ have the same score list if and only if they are $\{C_3, C_4\}$-equivalent.
\end{cor}
\noindent\textit{Proof}: (i) This follows from Theorem \ref{thmA9.3.4}(i) since $\mathcal{L}_G\subseteq\{3\}$.
\indent\par (ii) $(\Rightarrow)$ This follows from Theorem \ref{thmA9.3.1}, (i) and the Refinement Lemma.
\qed

\indent\par We remark that Corollary \ref{corA9.3.5} holds for $n$-partite tournaments since complete graphs are chordal graphs. Note via Examples \ref{egA9.3.6} and \ref{egA9.3.7} that the necessities of Theorems \ref{thmA9.1.1} and \ref{thmA9.1.2}, respectively, do not hold for $n$-partite tournaments, $n\ge 3$. So, Corollary \ref{corA9.3.5} balances the lack of such an analogue for $n$-partite tournaments by showing their $\{C_3, C_4\}$-equivalence and is best possible in the sense that neither of $C_3$ or $C_4$ can be omitted.

\begin{eg}\label{egA9.3.6}
To construct a pair of non-isomorphic $n$-partite tournaments, $n\ge 3$, with the same score list but not $C_3$-equivalent, we capitalise on the absence of odd cycles in bipartite tournaments. Let $D$ and $D'$ be two non-isomorphic bipartite tournaments of $K(V_1,V_2)$ with the same score list; thus they are $C_4$-equivalent by Theorem \ref{thmA9.1.2}. Let $F$ be the $n$-partite tournament such that $\langle V(F)-\bigcup\limits_{r=3}^n V_r\rangle_{F}\cong D$ and $V_i\rightarrow V_j$ whenever $i>j$ and $(i,j)\neq (2,1)$. Perform the same on $D'$ to obtain $F'$. Hence, $F$ and $F'$ have the same score list and are both $C_3$-free, i.e., they are not $C_3$-equivalent.
\end{eg}

\begin{eg} \label{egA9.3.7}
In Figure \ref{figA9.3.1}, $D$ and $D'$ are two non-isomorphic tripartite tournaments with the same score list. It is easy to see that $D$ is $C_4$-free, and thus $D$ and $D'$ are not $C_4$-equivalent. Similar to the previous example, $D$ and $D'$ can be extended to be examples of $n$-partite tournaments, $n\ge 4$.
\end{eg}
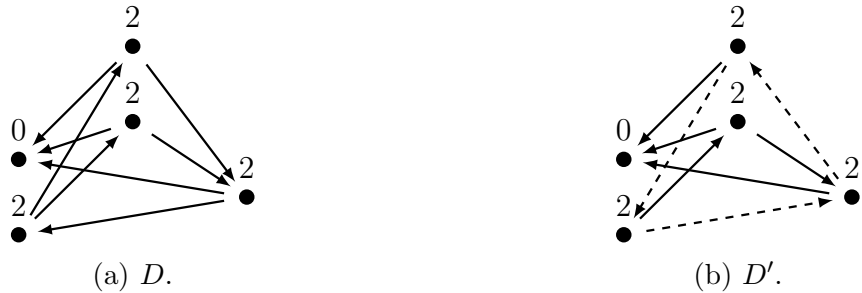
\begin{figure}[h]
\begin{subfigure}{.5\textwidth}
\begin{center}
\tikzstyle{every node}=[circle, draw, fill=black!100,
                       inner sep=0pt, minimum width=5pt]
\begin{tikzpicture}[thick,scale=0.5]%
\draw(0,0)node[label={[yshift=0.1cm]90:{$2$}}](12){};
\draw(0,2)node[label={[yshift=0.1cm]90:{$2$}}](11){};
\draw(-3,-1)node[label={[yshift=0.1cm]90:{$0$}}](31){};
\draw(-3,-3)node[label={[yshift=0.1cm]90:{$2$}}](32){};
\draw(3,-2)node[label={[yshift=0.1cm]90:{$2$}}](21){};

\draw[->, line width=0.3mm, >=latex, shorten <= 0.2cm, shorten >= 0.15cm](12)--(21);
\draw[->, line width=0.3mm, >=latex, shorten <= 0.2cm, shorten >= 0.15cm](11)--(21);

\draw[->, line width=0.3mm, >=latex, shorten <= 0.2cm, shorten >= 0.15cm](21)--(31);
\draw[->, line width=0.3mm, >=latex, shorten <= 0.2cm, shorten >= 0.15cm](21)--(32);

\draw[->, line width=0.3mm, >=latex, shorten <= 0.2cm, shorten >= 0.15cm](11)--(31);
\draw[ ->, line width=0.3mm, >=latex, shorten <= 0.2cm, shorten >= 0.15cm](12)--(31);
\draw[->, line width=0.3mm, >=latex, shorten <= 0.2cm, shorten >= 0.15cm](32)--(11);
\draw[->, line width=0.3mm, >=latex, shorten <= 0.2cm, shorten >= 0.15cm](32)--(12);
\end{tikzpicture}
{\caption{$D$.}}
\end{center}
\end{subfigure}%
\begin{subfigure}{0.5\textwidth}
\begin{center}
\tikzstyle{every node}=[circle, draw, fill=black!100,
                       inner sep=0pt, minimum width=5pt]
\begin{tikzpicture}[thick,scale=0.5]%
\draw(0,0)node[label={[yshift=0.1cm]90:{$2$}}](12){};
\draw(0,2)node[label={[yshift=0.1cm]90:{$2$}}](11){};
\draw(-3,-1)node[label={[yshift=0.1cm]90:{$0$}}](31){};
\draw(-3,-3)node[label={[yshift=0.1cm]90:{$2$}}](32){};
\draw(3,-2)node[label={[yshift=0.1cm]90:{$2$}}](21){};

\draw[dashed, ->, line width=0.3mm, >=latex, shorten <= 0.2cm, shorten >= 0.15cm](21)--(11);
\draw[->, line width=0.3mm, >=latex, shorten <= 0.2cm, shorten >= 0.15cm](12)--(21);
\draw[dashed, ->, line width=0.3mm, >=latex, shorten <= 0.2cm, shorten >= 0.15cm](32)--(21);
\draw[->, line width=0.3mm, >=latex, shorten <= 0.2cm, shorten >= 0.15cm](21)--(31);

\draw[->, line width=0.3mm, >=latex, shorten <= 0.2cm, shorten >= 0.15cm](32)--(12);
\draw[->, line width=0.3mm, >=latex, shorten <= 0.2cm, shorten >= 0.15cm](12)--(31);
\draw[->, line width=0.3mm, >=latex, shorten <= 0.2cm, shorten >= 0.15cm](11)--(31);
\draw[dashed, ->, line width=0.3mm, >=latex, shorten <= 0.2cm, shorten >= 0.15cm](11)--(32);
\end{tikzpicture}
{\caption{$D'$.}}
\end{center}
\end{subfigure}
{\caption{Two tripartite tournaments with the same score list but are not $C_4$-equivalent.}\label{figA9.3.1}}
\end{figure}

\section{Orientations with score-list parity}
In this section, we investigate orientations whose score lists preserve parity. To put our results in context, we state two results by Waldrop \cite{CW2}. Let $TT_3$ denote the transitive tournament of order $3$.
\begin{lem} (Waldrop \cite{CW2}) \label{lemA9.4.1}
In the family of tournaments with order at least $4$, $TT_3$ refines $C_3$.
\end{lem}

\begin{thm} (Waldrop \cite{CW2}) \label{thmA9.4.2}
Two tournaments of order $n$ have the same number of vertices of even score (equivalently, odd score) if and only if they are $TT_3$-equivalent.
\end{thm}

\noindent\par Let $G$ be a graph and $D$ and $D'$ be two orientations of $G(p_1, p_2, \ldots, p_n)$ with score lists $S$ and $S'$ respectively. We say that $D$ and $D'$ have \textit{score-list parity} if there exists some automorphism $f$ of $G(p_1, p_2, \ldots, p_n)$ such that $od_{D}(v)\equiv od_{D'}(f(v)) \pmod{2}$ for all $v\in V(D)$. In particular, $D$ and $D'$ have the same number of vertices with even score (equivalently, odd score).
\noindent\par Since a $TT_3$-reversal preserves score parity (and $TT_3$ is self-converse), the sufficiency of Theorem \ref{thmA9.4.2} is trivially true. As mentioned, $\mathscr{C}_k$ for $k\ge 3$ is self-converse and the next lemma shows that a $\mathscr{C}_k$-reversal preserves score parity. Therefore, most of the sufficiency proofs of the following results are trivial and will be omitted.

\begin{lem}\label{lemA9.4.3}
For $k\ge 3$, a $\mathscr{C}_k$-reversal preserves score parity.
\end{lem}
\noindent\textit{Proof}: Let $D\in \mathscr{C}_k$. Observe for any $v\in V(D)$, $od_D(v)=1$ if and only if $od_{\tilde{D}}(v)=1$, and $od_D(v)=2$ if and only if $od_{\tilde{D}}(v)=0$. Hence, the lemma follows.
\qed

\begin{thm}\label{thmA9.4.4}
Let $G$ be a graph. Two orientations of $G(p_1,p_2, \ldots, p_n)$ have score-list parity if and only if each orientation can be obtained from the other by successively reversing the arcs of orientations of pairwise edge-disjoint $k$-cycles, $k\ge 3$.
\end{thm}
\noindent\textit{Proof}: $(\Rightarrow)$ Let $D$ and $D'$ be two orientations of $G(p_1, p_2, \ldots, p_n)$ with score-list parity. By definition, there exists some automorphism $f$ of $G(p_1, p_2, \ldots, p_n)$ such that $od_{D}(v)\equiv od_{D'}(f(v)) \pmod{2}$ for all $v\in V(D)$. By Lemma \ref{lemA9.2.3} and denoting $\mathcal{D}=DDG(f;D,D')$, we have
\begin{align*}
od_\mathcal{D}(v)-id_\mathcal{D}(v)&=od_{D}(v)-od_{D'}(f(v))\equiv 0 \pmod{2} \text{ and thus},
\\deg_{\mathrm{U}(\mathcal{D})}(v)&=od_\mathcal{D}(v)+id_\mathcal{D}(v)\equiv 0 \pmod{2}
\end{align*}
for all $v\in V(\mathcal{D})$. That is, every vertex in $V(\mathrm{U}(\mathcal{D}))$ has even degree in $\mathrm{U}(\mathcal{D})$. So, $\mathrm{U}(\mathcal{D})$ is an Eulerian graph and can be partitioned into edge-disjoint cycles.
\indent\par $(\Leftarrow)$ This follows from Lemma \ref{lemA9.4.3}.
\qed

\begin{cor}\label{corA9.4.5}
Let $G$ be a graph with exactly one automorphism, namely the identity map. If every vertex has the same score parity in two orientations of $G$, then each orientation can be obtained from the other by successively reversing the arcs of orientations of pairwise edge-disjoint $k$-cycles, $k\ge 3$. Moreover, any bridge of $G$ is oriented the same in both orientations.
\end{cor}
\noindent\textit{Proof}: The first part follows by invoking Theorem \ref{thmA9.4.4} with $p_i=1$ for all $1\le i\le n$ and $f$ as the identity map of $V(G)$. The ``Moreover" part follows from the fact that an orientation of a bridge is invariant under any $\mathscr{C}_k$-reversal.
\qed

\indent\par We remark that if $G$ is a tree satisfying the conditions of Corollary \ref{corA9.4.5}, then the two said orientations are identical, not merely isomorphic. The next example shows that Corollary \ref{corA9.4.5} fails to hold if the two orientations merely have the same number of vertices of even score (equivalently, odd score).
\begin{eg}\label{egA9.4.6}
In Figure \ref{figA9.4.2}, $D$ and $D'$ are two non-isomorphic orientations of a graph with exactly one automorphism. They have the same number of vertices of even score and the bridge $uv$ oriented differently.
\end{eg}
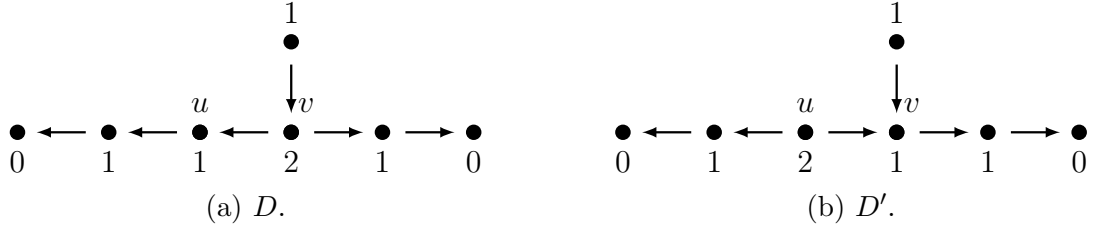
\begin{figure}[h]
\begin{subfigure}{.5\textwidth}
\begin{center}
\tikzstyle{every node}=[circle, draw, fill=black!100,
                       inner sep=0pt, minimum width=5pt]
\begin{tikzpicture}[thick,scale=0.6]%
\draw(0,0)node[label={[yshift=-0.1cm]270:{$0$}}](1){};
\draw(2,0)node[label={[yshift=-0.1cm]270:{$1$}}](2){};
\draw(4,0)node[label={[yshift=-0.1cm]270:{$1$}}](3){};
\draw(4,0)node[label={[yshift=0.1cm]90:{$u$}}](3){};
\draw(6,0)node[label={[yshift=-0.1cm]270:{$2$}}](4){};
\draw(6,0)node[label={[xshift=0.2cm, yshift=0.1cm]90:{$v$}}](4){};
\draw(8,0)node[label={[yshift=-0.1cm]270:{$1$}}](5){};
\draw(10,0)node[label={[yshift=-0.1cm]270:{$0$}}](6){};
\draw(6,2)node[label={[yshift=0.1cm]90:{$1$}}](7){};

\draw[->, line width=0.3mm, >=latex, shorten <= 0.2cm, shorten >= 0.15cm](4)--(3);
\draw[->, line width=0.3mm, >=latex, shorten <= 0.2cm, shorten >= 0.15cm](3)--(2);
\draw[->, line width=0.3mm, >=latex, shorten <= 0.2cm, shorten >= 0.15cm](2)--(1);
\draw[->, line width=0.3mm, >=latex, shorten <= 0.2cm, shorten >= 0.15cm](4)--(5);
\draw[->, line width=0.3mm, >=latex, shorten <= 0.2cm, shorten >= 0.15cm](5)--(6);
\draw[->, line width=0.3mm, >=latex, shorten <= 0.2cm, shorten >= 0.15cm](7)--(4);
\end{tikzpicture}
{\caption{$D$.}}
\end{center}
\end{subfigure}%
\begin{subfigure}{0.5\textwidth}
\begin{center}
\tikzstyle{every node}=[circle, draw, fill=black!100,
                       inner sep=0pt, minimum width=5pt]
\begin{tikzpicture}[thick,scale=0.6]%
\draw(0,0)node[label={[yshift=-0.1cm]270:{$0$}}](1){};
\draw(2,0)node[label={[yshift=-0.1cm]270:{$1$}}](2){};
\draw(4,0)node[label={[yshift=-0.1cm]270:{$2$}}](3){};
\draw(4,0)node[label={[yshift=0.1cm]90:{$u$}}](3){};
\draw(6,0)node[label={[yshift=-0.1cm]270:{$1$}}](4){};
\draw(6,0)node[label={[xshift=0.2cm, yshift=0.1cm]90:{$v$}}](4){};
\draw(8,0)node[label={[yshift=-0.1cm]270:{$1$}}](5){};
\draw(10,0)node[label={[yshift=-0.1cm]270:{$0$}}](6){};
\draw(6,2)node[label={[yshift=0.1cm]90:{$1$}}](7){};

\draw[->, line width=0.3mm, >=latex, shorten <= 0.2cm, shorten >= 0.15cm](3)--(4);
\draw[->, line width=0.3mm, >=latex, shorten <= 0.2cm, shorten >= 0.15cm](3)--(2);
\draw[->, line width=0.3mm, >=latex, shorten <= 0.2cm, shorten >= 0.15cm](2)--(1);
\draw[->, line width=0.3mm, >=latex, shorten <= 0.2cm, shorten >= 0.15cm](4)--(5);
\draw[->, line width=0.3mm, >=latex, shorten <= 0.2cm, shorten >= 0.15cm](5)--(6);
\draw[->, line width=0.3mm, >=latex, shorten <= 0.2cm, shorten >= 0.15cm](7)--(4);
\end{tikzpicture}
{\caption{$D'$.}}
\end{center}
\end{subfigure}
{\caption{Two orientations with the same number of vertices of even score but are not $\bigcup\limits_{k\ge 3}\mathscr{C}_k$-equivalent.}\label{figA9.4.2}}
\end{figure}

\noindent\par As in Section 3, we aim to strengthen Theorem \ref{thmA9.4.4} by dropping the edge-disjoint condition. We start by examining orientations of tree vertex-multiplications which encompass the bipartite tournaments.

\begin{ppn}\label{ppnA9.4.7} Let $T$ be a tree.
\\(i) In the family of orientations of $T(p_1, p_2, \ldots, p_n)$, $\mathscr{C}_4$ refines $\bigcup\limits_{k\ge 2}\mathscr{C}_{2k}$.
\\(ii) Two orientations of $T(p_1, p_2, \ldots, p_n)$ have score-list parity if and only if they are $\mathscr{C}_4$-equivalent.
\end{ppn}
\noindent\textit{Proof}: (i) Let $D$ be an orientation of $T(p_1, p_2, \ldots, p_n)$ containing an orientation $F$ of a $2k$-cycle, say $u_1 u_2 \ldots u_{2k} u_1$, $k\ge 2$. Recall that $T(p_1, p_2, \ldots, p_n)$ is bipartite and thus has no odd cycles. If $k=2$, then we are done. So, consider $k\ge 3$. $T$ is a tree implies that $F$ contains some pair of vertices $u_i$ and $u_j$ from the same partite set. Hence by relabelling the vertices if necessary, we may assume $1\le i<j\le 2k$ and invoke the Generalised Reduction$^*$ Lemma (repeatedly) to deduce (i).
\indent\par (ii) $(\Rightarrow)$ By Theorem \ref{thmA9.4.4}, two orientations of $T(p_1, p_2, \ldots, p_n)$ with score-list parity are $\bigcup\limits_{k\ge 3}\mathscr{C}_k$-equivalent. So, they are $\mathscr{C}_4$-equivalent by (i), the absence of odd cycles in $T(p_1, p_2, \ldots, p_n)$ and the Refinement Lemma.
\qed

\indent\par Parallel to Section 3, the following generalisation of Proposition \ref{ppnA9.4.7} can be derived. Its proof is very much alike to Theorem \ref{thmA9.3.4} and therefore omitted.

\begin{thm}\label{thmA9.4.8}
Let $G$ be a graph.
\\(i) In the family of orientations of $G(p_1, p_2, \ldots, p_n)$, $\bigcup\limits_{i\in \{4\}\cup \mathcal{L}_G}\mathscr{C}_i$ refines $\bigcup\limits_{k\ge 3}\mathscr{C}_k$.
\\(ii) Two orientations of $G(p_1,p_2, \ldots, p_n)$ have score-list parity if and only if they are $\bigcup\limits_{i\in \{4\}\cup \mathcal{L}_G}\mathscr{C}_i$-equivalent.
\end{thm}
{\color{blue} Proof omitted from journal paper
\\\noindent\textit{Proof}: (i) Let $D$ be an orientation of $G(p_1, p_2, \ldots, p_n)$ containing an orientation $F$ of a $k$-cycle $u_1 u_2\ldots u_k u_1$, $k\ge 3$. We start by proving a claim.
\\
\\Claim: $\bigcup\limits_{i\in \{4\}\cup \mathcal{L}^*_G}\mathscr{C}_i$ refines $\bigcup\limits_{k\ge 3}\mathscr{C}_k$.
\indent\par If $k=3$, then $\mathrm{UPG}(F)$ is a $3$-cycle in $G$, i.e., $k\in \mathcal{L}^*_G$. As before, it suffices to consider $k\ge 5$. We proceed by induction on the number of pairs $p$ of vertices $u_i, u_j$ of $F$ in a common partite set. If $p=0$, then $\mathrm{UPG}(F)$ is a cycle in $G$ and $k\in\mathcal{L}^*_G$. Suppose $p\ge 1$. By relabelling vertices if necessary, we may assume $u_i, u_j\in V(F)\cap V_w$ for some $1\le i<j\le k$ and $w\in V(G)$; we further assume the difference $j-i$ to be the smallest among all such pairs. If $F^*=\langle \{u_i, u_{i+1},\ldots, u_j\}\rangle_F$ and $H=\mathrm{UPG}(F^*)$ is a tree in $G$, then by the Generalised Reduction$^*$ Lemma with $F'$ as given, an $F$-reversal is equivalent to some sequence of $\mathscr{C}_4$-reversals and an $F'$-reversal. Since $\mathrm{U}(F')$ is a cycle with at most $p-1$ pairs of vertices in a common partite set, the claim follows by induction hypothesis.
\noindent\par Hence, consider the case where $H$ contains some cycle in $G$. By the Generalised Reduction Lemma with $F'$ and $F''$ as labelled, an $F$-reversal is equivalent to an $F'$-reversal followed by an $F''$-reversal. Let the arc joining $u_i$ and $u_{j-1}$ in $D$ be $a_{i,j-1}$. By invoking the Generalised Reduction$^*$ Lemma on $F''$ again, we deduce an $F''$-reversal is equivalent to reversing the arcs of the orientations $F'''=\langle\{u_i, u_{i+1}, \ldots, u_{j-1}\}\rangle_{F''}\cup\{a_{i,j-1}\}$ followed by $\langle\{u_{i-1}, u_i, u_{j-1}, u_j\}\rangle_{F''}\cup\{\tilde{a}_{i,j-1}\}$ of cycles; note that the latter is a $\mathscr{C}_4$-reversal. By minimality of $j-i$, $\mathrm{UPG}(F''')$ is a cycle in $G$; and since $F'$ has at most $p-1$ pairs of vertices in a common partite set, the claim follows by induction hypothesis. This completes the proof of the claim.
\\
\indent\par By the above proof, we may assume that $\mathrm{UPG}(F)$ is a cycle in $G$ and no two vertices of $F$ belong to the same partite set. It remains to show that an $F$-reversal is equivalent to some sequence of $\bigcup\limits_{i\in \{4\}\cup \mathcal{L}_G}\mathscr{C}_i$-reversals. We proceed by induction on the number of chords in $\langle \mathrm{UP}(F)\rangle_G$. If $k= 4$ or $\langle \mathrm{UP}(F)\rangle_G$ is chordless (which includes the case $k=3$), then $k\in \{4\}\cup\mathcal{L}_G$ and we are done. So, consider $k\ge 5$ and $\langle \mathrm{UP}(F)\rangle_G$ has $t$ chords, $t\ge 1$. Let $a_{1,i}$ denote the arc joining $u_1$ and $u_i$ in $D$. WLOG, assume that there exists some $3\le i\le k-1$ such that $\mathrm{UPG}(a_{1,i})$ is a chord of $\langle \mathrm{UP}(F)\rangle_G$. We may further assume that $F_0=\langle \{u_i, u_{i+1},\ldots, u_k, u_1\}\rangle_F$ and $\langle \mathrm{UP}(F_0)\rangle_G$ is a chordless cycle in $G$. Then an $F$-reversal is equivalent to reversing the arcs of the orientations $F_1=\langle\{u_1, u_2, \ldots, u_i\}\rangle_F\cup \{a_{1,i}\}$ and $F_2=F_0 \cup \{\tilde{a}_{1,i}\}$ of cycles. Note that no two vertices of $F_2$ belong to the same partite set and $\langle \mathrm{UP}(F_2)\rangle_G$ is a chordless cycle in $G$, i.e., $F_2\in \bigcup\limits_{i\in \mathcal{L}_G}\mathscr{C}_i$. Since $\langle \mathrm{UP}(F_1)\rangle_G$ has at most $t-1$ chords, the result follows by induction hypothesis.
\indent\par (ii) $(\Rightarrow)$ This follows from Theorem \ref{thmA9.4.4}, (i) and the Refinement Lemma.
\qed
}

\begin{cor}\label{corA9.4.9}
Let $G$ be a chordal graph.
\\(i) In the family of orientations of $G(p_1, p_2, \ldots, p_n)$, $\mathscr{C}_3\cup \mathscr{C}_4$ refines $\bigcup\limits_{k\ge 3}\mathscr{C}_k$.
\\(ii) Two orientations of $G(p_1,p_2, \ldots, p_n)$ have score-list parity if and only if they are $\mathscr{C}_3\cup \mathscr{C}_4$-equivalent.
\end{cor}
\noindent\textit{Proof}: (i) This follows from Theorem \ref{thmA9.4.8}(i) since $\mathcal{L}_G\subseteq\{3\}$.
\indent\par (ii) $(\Rightarrow)$ This follows from Theorem \ref{thmA9.4.4}, (i) and the Refinement Lemma.
\qed

\indent\par Akin to Section 3, Corollary \ref{corA9.4.9} holds for $n$-partite tournaments. We show that this can be further strengthened by considering the cases $n\ge 4$ and $n=3$ separately. The next lemma caters to the case $n\ge 4$ (we show that it also holds for $n=3$).
\begin{lem}\label{lemA9.4.10}
In the family of $n$-partite tournaments with $n\ge 3$, $TT_3$ refines $C_4$.
\end{lem}
\noindent\textit{Proof}: Let $T$ be an $n$-partite tournament, $n\ge 3$, containing a $4$-dicycle $Z=u_1 u_2 u_3 u_4 u_1$ to be reversed. Suppose there exists some $i=1,2$, such that the vertices $u_i$ and $u_{i+2}$ are in distinct partite sets, say $i=1$. WLOG, we may assume $u_1\rightarrow u_3$. Then a $Z$-reversal is equivalent to the sequence $u_1\rightarrow u_2\rightarrow u_3 \leftarrow u_1$, $u_3\rightarrow u_4\rightarrow u_1 \leftarrow u_3$ of $TT_3$-reversals.
\noindent\par Next, assume $u_1$ and $u_3$ are in a partite set, say $V_1$, while $u_2$ and $u_4$ are in another, say $V_2$. Let $v\in V_3$. By symmetry and duality, it suffices to consider the following four cases. If $v\rightarrow \{u_1, u_2, u_3, u_4\}$, then a $Z$-reversal is equivalent to the sequence
\begin{align*}
v\rightarrow u_2\rightarrow u_3 \leftarrow v, v\rightarrow u_4 \rightarrow u_1\leftarrow v, u_1\rightarrow u_2\rightarrow v\leftarrow u_1, u_3\rightarrow u_4\rightarrow v\leftarrow u_3
\end{align*}
of $TT_3$-reversals (see Figure \ref{figA9.4.3}). If $\{u_1, u_3\}\rightarrow v\rightarrow\{u_2, u_4\}$, then a $Z$-reversal is equivalent to the sequence
\begin{align*}
u_3\rightarrow v\rightarrow u_4\leftarrow u_3, u_4\rightarrow u_1\rightarrow v\leftarrow u_4, v\rightarrow u_1\rightarrow u_2\leftarrow v, u_2\rightarrow v\rightarrow u_3\leftarrow u_2
\end{align*}
of $TT_3$-reversals (see Figure \ref{figA9.4.4}). If $u_4\rightarrow v\rightarrow\{u_1, u_2, u_3\}$, then a $Z$-reversal is equivalent to the sequence
\begin{align*}
v\rightarrow u_2\rightarrow u_3\leftarrow v, u_3\rightarrow u_4\rightarrow v\leftarrow u_3, v\rightarrow u_4\rightarrow u_1\leftarrow v, u_1\rightarrow u_2\rightarrow v\leftarrow u_1
\end{align*}
of $TT_3$-reversals (see Figure \ref{figA9.4.5}). If $\{u_1, u_2\}\rightarrow v\rightarrow\{u_3, u_4\}$, then a $Z$-reversal is equivalent to the sequence
\begin{align*}
v\rightarrow u_3\rightarrow u_4\leftarrow v, u_4\rightarrow u_1\rightarrow v\leftarrow u_4, u_2\rightarrow u_3\rightarrow v\leftarrow u_3, v\rightarrow u_1\rightarrow u_2\leftarrow v
\end{align*}
of $TT_3$-reversals (see Figure \ref{figA9.4.6}). This completes the proof.
\qed

\begin{center}
\tikzstyle{every node}=[circle, draw, fill=black!100,
                       inner sep=0pt, minimum width=6pt]
\begin{tikzpicture}[thick,scale=0.7]%
\draw(-5,2.5)node[label={[yshift=0cm] 90:{$v$}}](1v){};
\draw(-5,1)node[label={[yshift=0cm] 270:{$u_1$}}](1u1){};
\draw(-6,0)node[label={[yshift=0cm] 270:{$u_2$}}](1u2){};
\draw(-4,0)node[label={[yshift=0cm] 270:{$u_4$}}](1u4){};
\draw(-5,-1)node[label={[yshift=0cm] 270:{$u_3$}}](1u3){};

\draw(0,2.5)node[label={[yshift=0cm] 90:{$v$}}](2v){};
\draw(0,1)node[label={[yshift=0cm] 270:{$u_1$}}](2u1){};
\draw(-1,0)node[label={[yshift=0cm] 270:{$u_2$}}](2u2){};
\draw(1,0)node[label={[yshift=0cm] 270:{$u_4$}}](2u4){};
\draw(0,-1)node[label={[yshift=0cm] 270:{$u_3$}}](2u3){};

\draw(5,2.5)node[label={[yshift=0cm] 90:{$v$}}](3v){};
\draw(5,1)node[label={[yshift=0cm] 270:{$u_1$}}](3u1){};
\draw(4,0)node[label={[yshift=0cm] 270:{$u_2$}}](3u2){};
\draw(6,0)node[label={[yshift=0cm] 270:{$u_4$}}](3u4){};
\draw(5,-1)node[label={[yshift=0cm] 270:{$u_3$}}](3u3){};

\draw(10,2.5)node[label={[yshift=0cm] 90:{$v$}}](4v){};
\draw(10,1)node[label={[yshift=0cm] 270:{$u_1$}}](4u1){};
\draw(9,0)node[label={[yshift=0cm] 270:{$u_2$}}](4u2){};
\draw(11,0)node[label={[yshift=0cm] 270:{$u_4$}}](4u4){};
\draw(10,-1)node[label={[yshift=0cm] 270:{$u_3$}}](4u3){};

\draw(15,2.5)node[label={[yshift=0cm] 90:{$v$}}](5v){};
\draw(15,1)node[label={[yshift=0cm] 270:{$u_1$}}](5u1){};
\draw(14,0)node[label={[yshift=0cm] 270:{$u_2$}}](5u2){};
\draw(16,0)node[label={[yshift=0cm] 270:{$u_4$}}](5u4){};
\draw(15,-1)node[label={[yshift=0cm] 270:{$u_3$}}](5u3){};

\draw[->, line width=0.3mm, >=latex, shorten <= 0.1cm, shorten >= 0.05cm](1u1)--(1u2);
\draw[dashdotted, ->, line width=0.3mm, >=latex, shorten <= 0.1cm, shorten >= 0.05cm](1u2)--(1u3);
\draw[->, line width=0.3mm, >=latex, shorten <= 0.1cm, shorten >= 0.05cm](1u3)--(1u4);
\draw[->, line width=0.3mm, >=latex, shorten <= 0.1cm, shorten >= 0.05cm](1u4)--(1u1);
\draw[->, line width=0.3mm, >=latex, shorten <= 0.1cm, shorten >= 0.05cm](1v)--(1u1);
\draw[dashdotted, ->, line width=0.3mm, >=latex, shorten <= 0.1cm, shorten >= 0.05cm](1v)--(1u2);
\draw[dashdotted, ->, line width=0.3mm, >=latex, shorten <= 0.1cm, shorten >= 0.05cm](1v) to (-6.5, 2.5) to (-6.5,-1) to (1u3);
\draw[->, line width=0.3mm, >=latex, shorten <= 0.1cm, shorten >= 0.05cm](1v)--(1u4);

\draw[->, line width=0.3mm, >=latex, shorten <= 0.1cm, shorten >= 0.05cm](2u1)--(2u2);
\draw[->, line width=0.3mm, >=latex, shorten <= 0.1cm, shorten >= 0.05cm](2u3)--(2u2);
\draw[->, line width=0.3mm, >=latex, shorten <= 0.1cm, shorten >= 0.05cm](2u3)--(2u4);
\draw[dashdotted, ->, line width=0.3mm, >=latex, shorten <= 0.1cm, shorten >= 0.05cm](2u4)--(2u1);
\draw[dashdotted, ->, line width=0.3mm, >=latex, shorten <= 0.1cm, shorten >= 0.05cm](2v)--(2u1);
\draw[->, line width=0.3mm, >=latex, shorten <= 0.1cm, shorten >= 0.05cm](2u2)--(2v);
\draw[->, line width=0.3mm, >=latex, shorten <= 0.1cm, shorten >= 0.05cm](2u3) to (-1.5, -1) to (-1.5,2.5) to (2v);
\draw[dashdotted, ->, line width=0.3mm, >=latex, shorten <= 0.1cm, shorten >= 0.05cm](2v)--(2u4);

\draw[dashdotted, ->, line width=0.3mm, >=latex, shorten <= 0.1cm, shorten >= 0.05cm](3u1)--(3u2);
\draw[->, line width=0.3mm, >=latex, shorten <= 0.1cm, shorten >= 0.05cm](3u3)--(3u2);
\draw[->, line width=0.3mm, >=latex, shorten <= 0.1cm, shorten >= 0.05cm](3u3)--(3u4);
\draw[->, line width=0.3mm, >=latex, shorten <= 0.1cm, shorten >= 0.05cm](3u1)--(3u4);
\draw[dashdotted, ->, line width=0.3mm, >=latex, shorten <= 0.1cm, shorten >= 0.05cm](3u1)--(3v);
\draw[dashdotted, ->, line width=0.3mm, >=latex, shorten <= 0.1cm, shorten >= 0.05cm](3u2)--(3v);
\draw[->, line width=0.3mm, >=latex, shorten <= 0.1cm, shorten >= 0.05cm](3u3) to (3.5, -1) to (3.5,2.5) to (3v);
\draw[->, line width=0.3mm, >=latex, shorten <= 0.1cm, shorten >= 0.05cm](3u4)--(3v);

\draw[->, line width=0.3mm, >=latex, shorten <= 0.1cm, shorten >= 0.05cm](4u2)--(4u1);
\draw[->, line width=0.3mm, >=latex, shorten <= 0.1cm, shorten >= 0.05cm](4u3)--(4u2);
\draw[dashdotted, ->, line width=0.3mm, >=latex, shorten <= 0.1cm, shorten >= 0.05cm](4u3)--(4u4);
\draw[->, line width=0.3mm, >=latex, shorten <= 0.1cm, shorten >= 0.05cm](4u1)--(4u4);
\draw[->, line width=0.3mm, >=latex, shorten <= 0.1cm, shorten >= 0.05cm](4v)--(4u1);
\draw[->, line width=0.3mm, >=latex, shorten <= 0.1cm, shorten >= 0.05cm](4v)--(4u2);
\draw[dashdotted, ->, line width=0.3mm, >=latex, shorten <= 0.1cm, shorten >= 0.05cm](4u3) to (8.5, -1) to (8.5,2.5) to (4v);
\draw[dashdotted, ->, line width=0.3mm, >=latex, shorten <= 0.1cm, shorten >= 0.05cm](4u4)--(4v);

\draw[->, line width=0.3mm, >=latex, shorten <= 0.1cm, shorten >= 0.05cm](5u2)--(5u1);
\draw[->, line width=0.3mm, >=latex, shorten <= 0.1cm, shorten >= 0.05cm](5u3)--(5u2);
\draw[->, line width=0.3mm, >=latex, shorten <= 0.1cm, shorten >= 0.05cm](5u4)--(5u3);
\draw[->, line width=0.3mm, >=latex, shorten <= 0.1cm, shorten >= 0.05cm](5u1)--(5u4);
\draw[->, line width=0.3mm, >=latex, shorten <= 0.1cm, shorten >= 0.05cm](5v)--(5u1);
\draw[->, line width=0.3mm, >=latex, shorten <= 0.1cm, shorten >= 0.05cm](5v)--(5u2);
\draw[->, line width=0.3mm, >=latex, shorten <= 0.1cm, shorten >= 0.05cm](5v) to (13.5, 2.5) to (13.5,-1) to (5u3);
\draw[->, line width=0.3mm, >=latex, shorten <= 0.1cm, shorten >= 0.05cm](5v)--(5u4);

\draw [-angle 90, shorten >= 0.1cm, line join=round,decorate, decoration={zigzag, segment length=8,amplitude=1.8,post=lineto,post length=8pt}](-3,0.75)--(-2,0.75);
\draw [-angle 90, shorten >= 0.1cm, line join=round,decorate, decoration={zigzag, segment length=8,amplitude=1.8,post=lineto,post length=8pt}](2,0.75)--(3,0.75);
\draw [-angle 90, shorten >= 0.1cm, line join=round,decorate, decoration={zigzag, segment length=8,amplitude=1.8,post=lineto,post length=8pt}](7,0.75)--(8,0.75);
\draw [-angle 90, shorten >= 0.1cm, line join=round,decorate, decoration={zigzag, segment length=8,amplitude=1.8,post=lineto,post length=8pt}](12,0.75)--(13,0.75);
\end{tikzpicture}
{
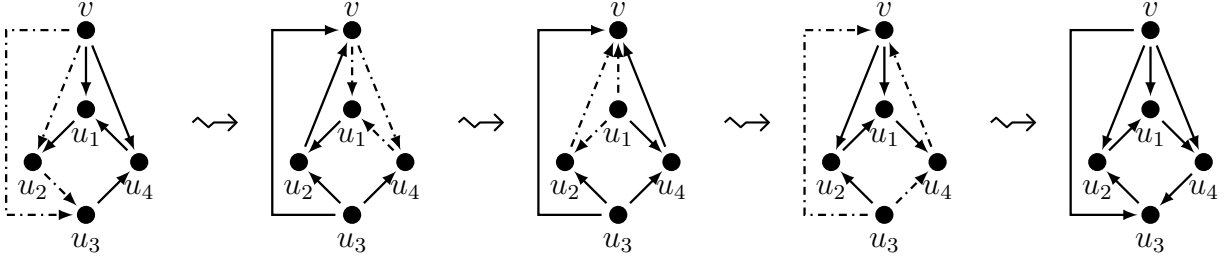
\captionof{figure}{$C_4$-reversal feasible via $TT_3$-reversals if $v\rightarrow \{u_1, u_2, u_3, u_4\}$.}\label{figA9.4.3}}

\tikzstyle{every node}=[circle, draw, fill=black!100,
                       inner sep=0pt, minimum width=6pt]
\begin{tikzpicture}[thick,scale=0.7]%
\draw(-5,2.5)node[label={[yshift=0cm] 90:{$v$}}](1v){};
\draw(-5,1)node[label={[yshift=0cm] 270:{$u_1$}}](1u1){};
\draw(-6,0)node[label={[yshift=0cm] 270:{$u_2$}}](1u2){};
\draw(-4,0)node[label={[yshift=0cm] 270:{$u_4$}}](1u4){};
\draw(-5,-1)node[label={[yshift=0cm] 270:{$u_3$}}](1u3){};

\draw(0,2.5)node[label={[yshift=0cm] 90:{$v$}}](2v){};
\draw(0,1)node[label={[yshift=0cm] 270:{$u_1$}}](2u1){};
\draw(-1,0)node[label={[yshift=0cm] 270:{$u_2$}}](2u2){};
\draw(1,0)node[label={[yshift=0cm] 270:{$u_4$}}](2u4){};
\draw(0,-1)node[label={[yshift=0cm] 270:{$u_3$}}](2u3){};

\draw(5,2.5)node[label={[yshift=0cm] 90:{$v$}}](3v){};
\draw(5,1)node[label={[yshift=0cm] 270:{$u_1$}}](3u1){};
\draw(4,0)node[label={[yshift=0cm] 270:{$u_2$}}](3u2){};
\draw(6,0)node[label={[yshift=0cm] 270:{$u_4$}}](3u4){};
\draw(5,-1)node[label={[yshift=0cm] 270:{$u_3$}}](3u3){};

\draw(10,2.5)node[label={[yshift=0cm] 90:{$v$}}](4v){};
\draw(10,1)node[label={[yshift=0cm] 270:{$u_1$}}](4u1){};
\draw(9,0)node[label={[yshift=0cm] 270:{$u_2$}}](4u2){};
\draw(11,0)node[label={[yshift=0cm] 270:{$u_4$}}](4u4){};
\draw(10,-1)node[label={[yshift=0cm] 270:{$u_3$}}](4u3){};

\draw(15,2.5)node[label={[yshift=0cm] 90:{$v$}}](5v){};
\draw(15,1)node[label={[yshift=0cm] 270:{$u_1$}}](5u1){};
\draw(14,0)node[label={[yshift=0cm] 270:{$u_2$}}](5u2){};
\draw(16,0)node[label={[yshift=0cm] 270:{$u_4$}}](5u4){};
\draw(15,-1)node[label={[yshift=0cm] 270:{$u_3$}}](5u3){};

\draw[->, line width=0.3mm, >=latex, shorten <= 0.1cm, shorten >= 0.05cm](1u1)--(1u2);
\draw[->, line width=0.3mm, >=latex, shorten <= 0.1cm, shorten >= 0.05cm](1u2)--(1u3);
\draw[dashdotted, ->, line width=0.3mm, >=latex, shorten <= 0.1cm, shorten >= 0.05cm](1u3)--(1u4);
\draw[->, line width=0.3mm, >=latex, shorten <= 0.1cm, shorten >= 0.05cm](1u4)--(1u1);
\draw[->, line width=0.3mm, >=latex, shorten <= 0.1cm, shorten >= 0.05cm](1u1)--(1v);
\draw[->, line width=0.3mm, >=latex, shorten <= 0.1cm, shorten >= 0.05cm](1v)--(1u2);
\draw[dashdotted, ->, line width=0.3mm, >=latex, shorten <= 0.1cm, shorten >= 0.05cm](1u3) to (-6.5, -1) to (-6.5,2.5) to (1v);
\draw[dashdotted, ->, line width=0.3mm, >=latex, shorten <= 0.1cm, shorten >= 0.05cm](1v)--(1u4);

\draw[->, line width=0.3mm, >=latex, shorten <= 0.1cm, shorten >= 0.05cm](2u1)--(2u2);
\draw[->, line width=0.3mm, >=latex, shorten <= 0.1cm, shorten >= 0.05cm](2u2)--(2u3);
\draw[->, line width=0.3mm, >=latex, shorten <= 0.1cm, shorten >= 0.05cm](2u4)--(2u3);
\draw[dashdotted, ->, line width=0.3mm, >=latex, shorten <= 0.1cm, shorten >= 0.05cm](2u4)--(2u1);
\draw[dashdotted, ->, line width=0.3mm, >=latex, shorten <= 0.1cm, shorten >= 0.05cm](2u1)--(2v);
\draw[->, line width=0.3mm, >=latex, shorten <= 0.1cm, shorten >= 0.05cm](2v)--(2u2);
\draw[->, line width=0.3mm, >=latex, shorten <= 0.1cm, shorten >= 0.05cm](2v) to (-1.5, 2.5) to (-1.5,-1) to (2u3);
\draw[dashdotted, ->, line width=0.3mm, >=latex, shorten <= 0.1cm, shorten >= 0.05cm](2u4)--(2v);

\draw[dashdotted, ->, line width=0.3mm, >=latex, shorten <= 0.1cm, shorten >= 0.05cm](3u1)--(3u2);
\draw[->, line width=0.3mm, >=latex, shorten <= 0.1cm, shorten >= 0.05cm](3u2)--(3u3);
\draw[->, line width=0.3mm, >=latex, shorten <= 0.1cm, shorten >= 0.05cm](3u4)--(3u3);
\draw[->, line width=0.3mm, >=latex, shorten <= 0.1cm, shorten >= 0.05cm](3u1)--(3u4);
\draw[dashdotted, ->, line width=0.3mm, >=latex, shorten <= 0.1cm, shorten >= 0.05cm](3v)--(3u1);
\draw[dashdotted, ->, line width=0.3mm, >=latex, shorten <= 0.1cm, shorten >= 0.05cm](3v)--(3u2);
\draw[->, line width=0.3mm, >=latex, shorten <= 0.1cm, shorten >= 0.05cm](3v) to (3.5, 2.5) to (3.5,-1) to (3u3);
\draw[->, line width=0.3mm, >=latex, shorten <= 0.1cm, shorten >= 0.05cm](3v)--(3u4);

\draw[->, line width=0.3mm, >=latex, shorten <= 0.1cm, shorten >= 0.05cm](4u2)--(4u1);
\draw[dashdotted, ->, line width=0.3mm, >=latex, shorten <= 0.1cm, shorten >= 0.05cm](4u2)--(4u3);
\draw[->, line width=0.3mm, >=latex, shorten <= 0.1cm, shorten >= 0.05cm](4u4)--(4u3);
\draw[->, line width=0.3mm, >=latex, shorten <= 0.1cm, shorten >= 0.05cm](4u1)--(4u4);
\draw[->, line width=0.3mm, >=latex, shorten <= 0.1cm, shorten >= 0.05cm](4u1)--(4v);
\draw[dashdotted, ->, line width=0.3mm, >=latex, shorten <= 0.1cm, shorten >= 0.05cm](4u2)--(4v);
\draw[dashdotted, ->, line width=0.3mm, >=latex, shorten <= 0.1cm, shorten >= 0.05cm](4v) to (8.5, 2.5) to (8.5,-1) to (4u3);
\draw[->, line width=0.3mm, >=latex, shorten <= 0.1cm, shorten >= 0.05cm](4v)--(4u4);

\draw[->, line width=0.3mm, >=latex, shorten <= 0.1cm, shorten >= 0.05cm](5u2)--(5u1);
\draw[->, line width=0.3mm, >=latex, shorten <= 0.1cm, shorten >= 0.05cm](5u3)--(5u2);
\draw[->, line width=0.3mm, >=latex, shorten <= 0.1cm, shorten >= 0.05cm](5u4)--(5u3);
\draw[->, line width=0.3mm, >=latex, shorten <= 0.1cm, shorten >= 0.05cm](5u1)--(5u4);
\draw[->, line width=0.3mm, >=latex, shorten <= 0.1cm, shorten >= 0.05cm](5u1)--(5v);
\draw[->, line width=0.3mm, >=latex, shorten <= 0.1cm, shorten >= 0.05cm](5v)--(5u2);
\draw[->, line width=0.3mm, >=latex, shorten <= 0.1cm, shorten >= 0.05cm](5u3) to (13.5, -1) to (13.5,2.5) to (5v);
\draw[->, line width=0.3mm, >=latex, shorten <= 0.1cm, shorten >= 0.05cm](5v)--(5u4);

\draw [-angle 90, shorten >= 0.1cm, line join=round,decorate, decoration={zigzag, segment length=8,amplitude=1.8,post=lineto,post length=8pt}](-3,0.75)--(-2,0.75);
\draw [-angle 90, shorten >= 0.1cm, line join=round,decorate, decoration={zigzag, segment length=8,amplitude=1.8,post=lineto,post length=8pt}](2,0.75)--(3,0.75);
\draw [-angle 90, shorten >= 0.1cm, line join=round,decorate, decoration={zigzag, segment length=8,amplitude=1.8,post=lineto,post length=8pt}](7,0.75)--(8,0.75);
\draw [-angle 90, shorten >= 0.1cm, line join=round,decorate, decoration={zigzag, segment length=8,amplitude=1.8,post=lineto,post length=8pt}](12,0.75)--(13,0.75);
\end{tikzpicture}
{
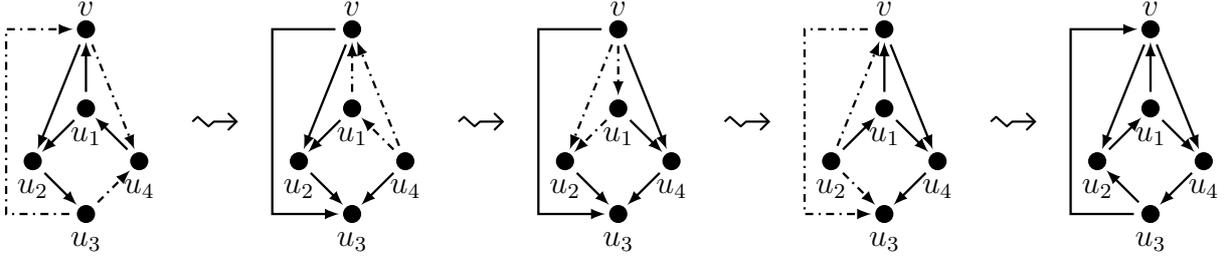
\captionof{figure}{$C_4$-reversal feasible via $TT_3$-reversals if $\{u_1, u_3\}\rightarrow v\rightarrow\{u_2, u_4\}$.}\label{figA9.4.4}}

\tikzstyle{every node}=[circle, draw, fill=black!100,
                       inner sep=0pt, minimum width=6pt]
\begin{tikzpicture}[thick,scale=0.7]%
\draw(-5,2.5)node[label={[yshift=0cm] 90:{$v$}}](1v){};
\draw(-5,1)node[label={[yshift=0cm] 270:{$u_1$}}](1u1){};
\draw(-6,0)node[label={[yshift=0cm] 270:{$u_2$}}](1u2){};
\draw(-4,0)node[label={[yshift=0cm] 270:{$u_4$}}](1u4){};
\draw(-5,-1)node[label={[yshift=0cm] 270:{$u_3$}}](1u3){};

\draw(0,2.5)node[label={[yshift=0cm] 90:{$v$}}](2v){};
\draw(0,1)node[label={[yshift=0cm] 270:{$u_1$}}](2u1){};
\draw(-1,0)node[label={[yshift=0cm] 270:{$u_2$}}](2u2){};
\draw(1,0)node[label={[yshift=0cm] 270:{$u_4$}}](2u4){};
\draw(0,-1)node[label={[yshift=0cm] 270:{$u_3$}}](2u3){};

\draw(5,2.5)node[label={[yshift=0cm] 90:{$v$}}](3v){};
\draw(5,1)node[label={[yshift=0cm] 270:{$u_1$}}](3u1){};
\draw(4,0)node[label={[yshift=0cm] 270:{$u_2$}}](3u2){};
\draw(6,0)node[label={[yshift=0cm] 270:{$u_4$}}](3u4){};
\draw(5,-1)node[label={[yshift=0cm] 270:{$u_3$}}](3u3){};

\draw(10,2.5)node[label={[yshift=0cm] 90:{$v$}}](4v){};
\draw(10,1)node[label={[yshift=0cm] 270:{$u_1$}}](4u1){};
\draw(9,0)node[label={[yshift=0cm] 270:{$u_2$}}](4u2){};
\draw(11,0)node[label={[yshift=0cm] 270:{$u_4$}}](4u4){};
\draw(10,-1)node[label={[yshift=0cm] 270:{$u_3$}}](4u3){};

\draw(15,2.5)node[label={[yshift=0cm] 90:{$v$}}](5v){};
\draw(15,1)node[label={[yshift=0cm] 270:{$u_1$}}](5u1){};
\draw(14,0)node[label={[yshift=0cm] 270:{$u_2$}}](5u2){};
\draw(16,0)node[label={[yshift=0cm] 270:{$u_4$}}](5u4){};
\draw(15,-1)node[label={[yshift=0cm] 270:{$u_3$}}](5u3){};

\draw[->, line width=0.3mm, >=latex, shorten <= 0.1cm, shorten >= 0.05cm](1u1)--(1u2);
\draw[dashdotted, ->, line width=0.3mm, >=latex, shorten <= 0.1cm, shorten >= 0.05cm](1u2)--(1u3);
\draw[->, line width=0.3mm, >=latex, shorten <= 0.1cm, shorten >= 0.05cm](1u3)--(1u4);
\draw[->, line width=0.3mm, >=latex, shorten <= 0.1cm, shorten >= 0.05cm](1u4)--(1u1);
\draw[->, line width=0.3mm, >=latex, shorten <= 0.1cm, shorten >= 0.05cm](1v)--(1u1);
\draw[dashdotted, ->, line width=0.3mm, >=latex, shorten <= 0.1cm, shorten >= 0.05cm](1v)--(1u2);
\draw[dashdotted, ->, line width=0.3mm, >=latex, shorten <= 0.1cm, shorten >= 0.05cm](1v) to (-6.5, 2.5) to (-6.5,-1) to (1u3);
\draw[->, line width=0.3mm, >=latex, shorten <= 0.1cm, shorten >= 0.05cm](1u4)--(1v);

\draw[->, line width=0.3mm, >=latex, shorten <= 0.1cm, shorten >= 0.05cm](2u1)--(2u2);
\draw[->, line width=0.3mm, >=latex, shorten <= 0.1cm, shorten >= 0.05cm](2u3)--(2u2);
\draw[dashdotted,->, line width=0.3mm, >=latex, shorten <= 0.1cm, shorten >= 0.05cm](2u3)--(2u4);
\draw[->, line width=0.3mm, >=latex, shorten <= 0.1cm, shorten >= 0.05cm](2u4)--(2u1);
\draw[->, line width=0.3mm, >=latex, shorten <= 0.1cm, shorten >= 0.05cm](2v)--(2u1);
\draw[->, line width=0.3mm, >=latex, shorten <= 0.1cm, shorten >= 0.05cm](2u2)--(2v);
\draw[dashdotted, ->, line width=0.3mm, >=latex, shorten <= 0.1cm, shorten >= 0.05cm](2u3) to (-1.5, -1) to (-1.5,2.5) to (2v);
\draw[dashdotted, ->, line width=0.3mm, >=latex, shorten <= 0.1cm, shorten >= 0.05cm](2u4)--(2v);

\draw[->, line width=0.3mm, >=latex, shorten <= 0.1cm, shorten >= 0.05cm](3u1)--(3u2);
\draw[->, line width=0.3mm, >=latex, shorten <= 0.1cm, shorten >= 0.05cm](3u3)--(3u2);
\draw[->, line width=0.3mm, >=latex, shorten <= 0.1cm, shorten >= 0.05cm](3u4)--(3u3);
\draw[dashdotted, ->, line width=0.3mm, >=latex, shorten <= 0.1cm, shorten >= 0.05cm](3u4)--(3u1);
\draw[dashdotted, ->, line width=0.3mm, >=latex, shorten <= 0.1cm, shorten >= 0.05cm](3v)--(3u1);
\draw[->, line width=0.3mm, >=latex, shorten <= 0.1cm, shorten >= 0.05cm](3u2)--(3v);
\draw[->, line width=0.3mm, >=latex, shorten <= 0.1cm, shorten >= 0.05cm](3v) to (3.5, 2.5) to (3.5, -1) to (3u3);
\draw[dashdotted, ->, line width=0.3mm, >=latex, shorten <= 0.1cm, shorten >= 0.05cm](3v)--(3u4);

\draw[dashdotted, ->, line width=0.3mm, >=latex, shorten <= 0.1cm, shorten >= 0.05cm](4u1)--(4u2);
\draw[->, line width=0.3mm, >=latex, shorten <= 0.1cm, shorten >= 0.05cm](4u3)--(4u2);
\draw[->, line width=0.3mm, >=latex, shorten <= 0.1cm, shorten >= 0.05cm](4u4)--(4u3);
\draw[->, line width=0.3mm, >=latex, shorten <= 0.1cm, shorten >= 0.05cm](4u1)--(4u4);
\draw[dashdotted, ->, line width=0.3mm, >=latex, shorten <= 0.1cm, shorten >= 0.05cm](4u1)--(4v);
\draw[dashdotted, ->, line width=0.3mm, >=latex, shorten <= 0.1cm, shorten >= 0.05cm](4u2)--(4v);
\draw[->, line width=0.3mm, >=latex, shorten <= 0.1cm, shorten >= 0.05cm](4v) to (8.5, 2.5) to (8.5, -1) to (4u3);
\draw[->, line width=0.3mm, >=latex, shorten <= 0.1cm, shorten >= 0.05cm](4u4)--(4v);

\draw[->, line width=0.3mm, >=latex, shorten <= 0.1cm, shorten >= 0.05cm](5u2)--(5u1);
\draw[->, line width=0.3mm, >=latex, shorten <= 0.1cm, shorten >= 0.05cm](5u3)--(5u2);
\draw[->, line width=0.3mm, >=latex, shorten <= 0.1cm, shorten >= 0.05cm](5u4)--(5u3);
\draw[->, line width=0.3mm, >=latex, shorten <= 0.1cm, shorten >= 0.05cm](5u1)--(5u4);
\draw[->, line width=0.3mm, >=latex, shorten <= 0.1cm, shorten >= 0.05cm](5v)--(5u1);
\draw[->, line width=0.3mm, >=latex, shorten <= 0.1cm, shorten >= 0.05cm](5v)--(5u2);
\draw[->, line width=0.3mm, >=latex, shorten <= 0.1cm, shorten >= 0.05cm](5v) to (13.5, 2.5) to (13.5, -1) to (5u3);
\draw[->, line width=0.3mm, >=latex, shorten <= 0.1cm, shorten >= 0.05cm](5u4)--(5v);

\draw [-angle 90, shorten >= 0.1cm, line join=round,decorate, decoration={zigzag, segment length=8,amplitude=1.8,post=lineto,post length=8pt}](-3,0.75)--(-2,0.75);
\draw [-angle 90, shorten >= 0.1cm, line join=round,decorate, decoration={zigzag, segment length=8,amplitude=1.8,post=lineto,post length=8pt}](2,0.75)--(3,0.75);
\draw [-angle 90, shorten >= 0.1cm, line join=round,decorate, decoration={zigzag, segment length=8,amplitude=1.8,post=lineto,post length=8pt}](7,0.75)--(8,0.75);
\draw [-angle 90, shorten >= 0.1cm, line join=round,decorate, decoration={zigzag, segment length=8,amplitude=1.8,post=lineto,post length=8pt}](12,0.75)--(13,0.75);
\end{tikzpicture}
{
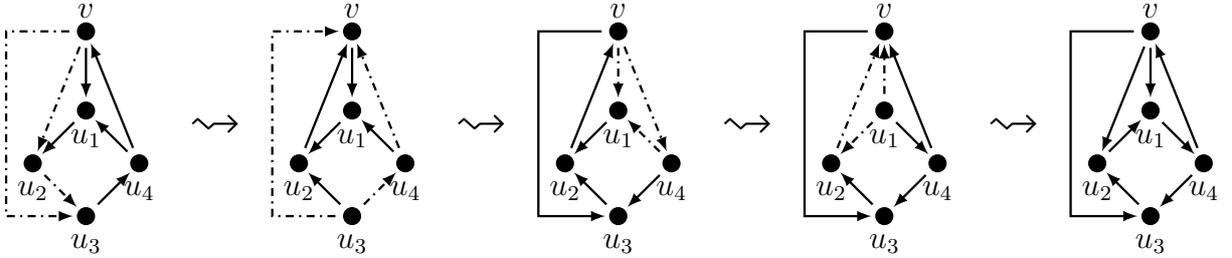
\captionof{figure}{$C_4$-reversal feasible via $TT_3$-reversals if $u_4\rightarrow v\rightarrow\{u_1, u_2, u_3\}$.}\label{figA9.4.5}}

\tikzstyle{every node}=[circle, draw, fill=black!100,
                       inner sep=0pt, minimum width=6pt]
\begin{tikzpicture}[thick,scale=0.7]%
\draw(-5,2.5)node[label={[yshift=0cm] 90:{$v$}}](1v){};
\draw(-5,1)node[label={[yshift=0cm] 270:{$u_1$}}](1u1){};
\draw(-6,0)node[label={[yshift=0cm] 270:{$u_2$}}](1u2){};
\draw(-4,0)node[label={[yshift=0cm] 270:{$u_4$}}](1u4){};
\draw(-5,-1)node[label={[yshift=0cm] 270:{$u_3$}}](1u3){};

\draw(0,2.5)node[label={[yshift=0cm] 90:{$v$}}](2v){};
\draw(0,1)node[label={[yshift=0cm] 270:{$u_1$}}](2u1){};
\draw(-1,0)node[label={[yshift=0cm] 270:{$u_2$}}](2u2){};
\draw(1,0)node[label={[yshift=0cm] 270:{$u_4$}}](2u4){};
\draw(0,-1)node[label={[yshift=0cm] 270:{$u_3$}}](2u3){};

\draw(5,2.5)node[label={[yshift=0cm] 90:{$v$}}](3v){};
\draw(5,1)node[label={[yshift=0cm] 270:{$u_1$}}](3u1){};
\draw(4,0)node[label={[yshift=0cm] 270:{$u_2$}}](3u2){};
\draw(6,0)node[label={[yshift=0cm] 270:{$u_4$}}](3u4){};
\draw(5,-1)node[label={[yshift=0cm] 270:{$u_3$}}](3u3){};

\draw(10,2.5)node[label={[yshift=0cm] 90:{$v$}}](4v){};
\draw(10,1)node[label={[yshift=0cm] 270:{$u_1$}}](4u1){};
\draw(9,0)node[label={[yshift=0cm] 270:{$u_2$}}](4u2){};
\draw(11,0)node[label={[yshift=0cm] 270:{$u_4$}}](4u4){};
\draw(10,-1)node[label={[yshift=0cm] 270:{$u_3$}}](4u3){};

\draw(15,2.5)node[label={[yshift=0cm] 90:{$v$}}](5v){};
\draw(15,1)node[label={[yshift=0cm] 270:{$u_1$}}](5u1){};
\draw(14,0)node[label={[yshift=0cm] 270:{$u_2$}}](5u2){};
\draw(16,0)node[label={[yshift=0cm] 270:{$u_4$}}](5u4){};
\draw(15,-1)node[label={[yshift=0cm] 270:{$u_3$}}](5u3){};

\draw[->, line width=0.3mm, >=latex, shorten <= 0.1cm, shorten >= 0.05cm](1u1)--(1u2);
\draw[->, line width=0.3mm, >=latex, shorten <= 0.1cm, shorten >= 0.05cm](1u2)--(1u3);
\draw[dashdotted, ->, line width=0.3mm, >=latex, shorten <= 0.1cm, shorten >= 0.05cm](1u3)--(1u4);
\draw[->, line width=0.3mm, >=latex, shorten <= 0.1cm, shorten >= 0.05cm](1u4)--(1u1);
\draw[->, line width=0.3mm, >=latex, shorten <= 0.1cm, shorten >= 0.05cm](1u1)--(1v);
\draw[->, line width=0.3mm, >=latex, shorten <= 0.1cm, shorten >= 0.05cm](1u2)--(1v);
\draw[dashdotted, ->, line width=0.3mm, >=latex, shorten <= 0.1cm, shorten >= 0.05cm](1v) to (-6.5, 2.5) to (-6.5,-1) to (1u3);
\draw[dashdotted, ->, line width=0.3mm, >=latex, shorten <= 0.1cm, shorten >= 0.05cm](1v)--(1u4);

\draw[->, line width=0.3mm, >=latex, shorten <= 0.1cm, shorten >= 0.05cm](2u1)--(2u2);
\draw[->, line width=0.3mm, >=latex, shorten <= 0.1cm, shorten >= 0.05cm](2u2)--(2u3);
\draw[->, line width=0.3mm, >=latex, shorten <= 0.1cm, shorten >= 0.05cm](2u4)--(2u3);
\draw[dashdotted, ->, line width=0.3mm, >=latex, shorten <= 0.1cm, shorten >= 0.05cm](2u4)--(2u1);
\draw[dashdotted, ->, line width=0.3mm, >=latex, shorten <= 0.1cm, shorten >= 0.05cm](2u1)--(2v);
\draw[->, line width=0.3mm, >=latex, shorten <= 0.1cm, shorten >= 0.05cm](2u2)--(2v);
\draw[->, line width=0.3mm, >=latex, shorten <= 0.1cm, shorten >= 0.05cm](2u3) to (-1.5, -1) to (-1.5, 2.5) to (2v);
\draw[dashdotted, ->, line width=0.3mm, >=latex, shorten <= 0.1cm, shorten >= 0.05cm](2u4)--(2v);

\draw[->, line width=0.3mm, >=latex, shorten <= 0.1cm, shorten >= 0.05cm](3u1)--(3u2);
\draw[dashdotted, ->, line width=0.3mm, >=latex, shorten <= 0.1cm, shorten >= 0.05cm](3u2)--(3u3);
\draw[->, line width=0.3mm, >=latex, shorten <= 0.1cm, shorten >= 0.05cm](3u4)--(3u3);
\draw[->, line width=0.3mm, >=latex, shorten <= 0.1cm, shorten >= 0.05cm](3u1)--(3u4);
\draw[->, line width=0.3mm, >=latex, shorten <= 0.1cm, shorten >= 0.05cm](3v)--(3u1);
\draw[dashdotted, ->, line width=0.3mm, >=latex, shorten <= 0.1cm, shorten >= 0.05cm](3u2)--(3v);
\draw[dashdotted, ->, line width=0.3mm, >=latex, shorten <= 0.1cm, shorten >= 0.05cm](3u3) to (3.5, -1) to (3.5, 2.5) to (3v);
\draw[->, line width=0.3mm, >=latex, shorten <= 0.1cm, shorten >= 0.05cm](3v)--(3u4);

\draw[dashdotted, ->, line width=0.3mm, >=latex, shorten <= 0.1cm, shorten >= 0.05cm](4u1)--(4u2);
\draw[->, line width=0.3mm, >=latex, shorten <= 0.1cm, shorten >= 0.05cm](4u3)--(4u2);
\draw[->, line width=0.3mm, >=latex, shorten <= 0.1cm, shorten >= 0.05cm](4u4)--(4u3);
\draw[->, line width=0.3mm, >=latex, shorten <= 0.1cm, shorten >= 0.05cm](4u1)--(4u4);
\draw[dashdotted, ->, line width=0.3mm, >=latex, shorten <= 0.1cm, shorten >= 0.05cm](4v)--(4u1);
\draw[dashdotted, ->, line width=0.3mm, >=latex, shorten <= 0.1cm, shorten >= 0.05cm](4v)--(4u2);
\draw[->, line width=0.3mm, >=latex, shorten <= 0.1cm, shorten >= 0.05cm](4v) to (8.5, 2.5) to (8.5, -1) to (4u3);
\draw[->, line width=0.3mm, >=latex, shorten <= 0.1cm, shorten >= 0.05cm](4v)--(4u4);

\draw[->, line width=0.3mm, >=latex, shorten <= 0.1cm, shorten >= 0.05cm](5u2)--(5u1);
\draw[->, line width=0.3mm, >=latex, shorten <= 0.1cm, shorten >= 0.05cm](5u3)--(5u2);
\draw[->, line width=0.3mm, >=latex, shorten <= 0.1cm, shorten >= 0.05cm](5u4)--(5u3);
\draw[->, line width=0.3mm, >=latex, shorten <= 0.1cm, shorten >= 0.05cm](5u1)--(5u4);
\draw[->, line width=0.3mm, >=latex, shorten <= 0.1cm, shorten >= 0.05cm](5u1)--(5v);
\draw[->, line width=0.3mm, >=latex, shorten <= 0.1cm, shorten >= 0.05cm](5u2)--(5v);
\draw[->, line width=0.3mm, >=latex, shorten <= 0.1cm, shorten >= 0.05cm](5v) to (13.5, 2.5) to (13.5, -1) to (5u3);
\draw[->, line width=0.3mm, >=latex, shorten <= 0.1cm, shorten >= 0.05cm](5v)--(5u4);

\draw [-angle 90, shorten >= 0.1cm, line join=round,decorate, decoration={zigzag, segment length=8,amplitude=1.8,post=lineto,post length=8pt}](-3,0.75)--(-2,0.75);
\draw [-angle 90, shorten >= 0.1cm, line join=round,decorate, decoration={zigzag, segment length=8,amplitude=1.8,post=lineto,post length=8pt}](2,0.75)--(3,0.75);
\draw [-angle 90, shorten >= 0.1cm, line join=round,decorate, decoration={zigzag, segment length=8,amplitude=1.8,post=lineto,post length=8pt}](7,0.75)--(8,0.75);
\draw [-angle 90, shorten >= 0.1cm, line join=round,decorate, decoration={zigzag, segment length=8,amplitude=1.8,post=lineto,post length=8pt}](12,0.75)--(13,0.75);
\end{tikzpicture}
{
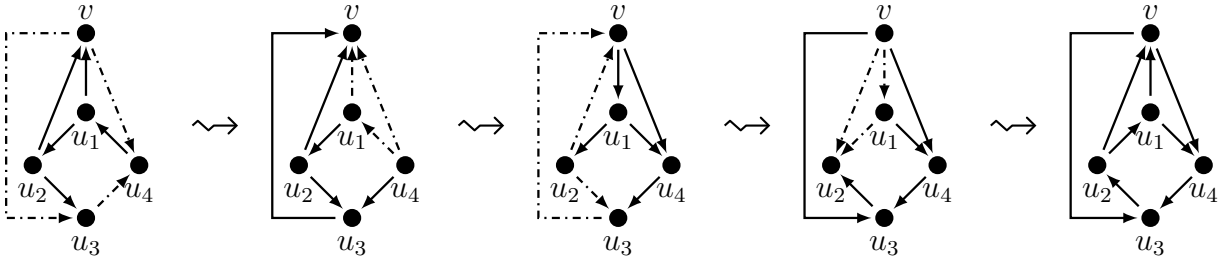
\captionof{figure}{$C_4$-reversal feasible via $TT_3$-reversals if $\{u_1, u_2\}\rightarrow v\rightarrow\{u_3, u_4\}$.}\label{figA9.4.6}}
\end{center}

\begin{thm}\label{thmA9.4.11} Let $n\ge 4$.
\\(i) In the family of $n$-partite tournaments, $TT_3$ refines $\bigcup\limits_{k\ge 3}\mathscr{C}_k$.
\\(ii) Two $n$-partite tournaments have score-list parity if and only if they are $TT_3$-equivalent.
\end{thm}
\noindent\textit{Proof}: (i) Let $T$ be an $n$-partite tournament, $n\ge 4$, containing an orientation $F$ of a $k$-cycle, say $u_1 u_2 \ldots u_k u_1$, $k\ge 3$. By Corollary \ref{corA9.4.9}(i), it suffices to consider $k=3$ or $4$, and show that an $F$-reversal in $T$ is equivalent to some sequence of $TT_3$-reversals. We proceed by induction on the number of arcs $t$ of $F$ that differ from the dicycle $C_k=u_1 u_2 \ldots u_k u_1$, i.e., $t=|\{i\mid u_{i+1} u_i\in A(F)\}|$ with addition taken modulo $k$. If $t=0$, then $F=C_k$. By Lemmas \ref{lemA9.4.1} and \ref {lemA9.4.10} respectively, any $C_3$-reversal and $C_4$-reversal are feasible in $T$ via $TT_3$-reversals.
\indent\par Suppose $t\ge 1$ and WLOG that $u_2 u_1\in A(F)-A(C_k)$. Since $n\ge 4$, there exists a vertex $v\in V(T)$ in a different partite set as $u_1$ and $u_2$. If $v\rightarrow u_2\rightarrow u_1\leftarrow v$ (or $u_2\rightarrow u_1\rightarrow v \leftarrow u_2$, or $u_2\rightarrow v\rightarrow u_1\leftarrow u_2$, or $u_2\rightarrow u_1\rightarrow v\rightarrow u_2$ resp.), then reversing the arcs of this $TT_3$ ($TT_3$, $TT_3$, $C_3$ resp.) yields an orientation $F^*$ of the $k$-cycle such that $|\{i\mid u_{i+1} u_i\in A(F^*)\cap A(C_k)\}|<t$. By induction hypothesis, an $F^*$-reversal in $T$ is feasible via $TT_3$-reversals. Next, reversing the arcs in $u_2\rightarrow u_1 \rightarrow v\leftarrow u_2$ ($v\rightarrow u_2\rightarrow u_1\leftarrow v$, $u_2\rightarrow u_1\rightarrow v\rightarrow u_2$, $u_2\rightarrow v\rightarrow u_1\leftarrow u_2$ resp.) completes exactly the $F$-reversal.
\indent\par (ii) $(\Rightarrow)$ This follows from Theorem \ref{thmA9.4.4}, (i) and the Refinement Lemma.
\qed

\indent\par By setting all partite sets to size one, Theorem \ref{thmA9.4.11}(ii) reduces to Theorem \ref{thmA9.4.2} for $n\ge 4$. The next example shows that the necessity of Theorem \ref{thmA9.4.11}(ii) does not hold generally for tripartite tournaments.
\begin{eg} \label{egA9.4.12}
In Figure \ref{figA9.4.7}, $D$ and $D'$ are two non-isomorphic tripartite tournaments with score-list parity. Furthermore, $D$ is $TT_3$-free.
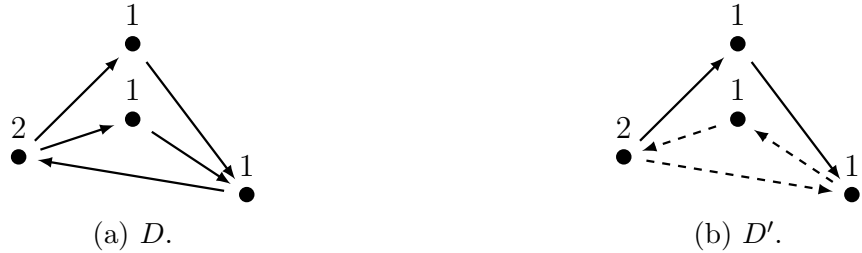
\begin{figure}[h]
\begin{subfigure}{.5\textwidth}
\begin{center}
\tikzstyle{every node}=[circle, draw, fill=black!100,
                       inner sep=0pt, minimum width=5pt]
\begin{tikzpicture}[thick,scale=0.5]%
\draw(0,0)node[label={[yshift=0.1cm]90:{$1$}}](12){};
\draw(0,2)node[label={[yshift=0.1cm]90:{$1$}}](11){};
\draw(-3,-1)node[label={[yshift=0.1cm]90:{$2$}}](31){};
\draw(3,-2)node[label={[yshift=0.1cm]90:{$1$}}](21){};

\draw[->, line width=0.3mm, >=latex, shorten <= 0.2cm, shorten >= 0.15cm](12)--(21);
\draw[->, line width=0.3mm, >=latex, shorten <= 0.2cm, shorten >= 0.15cm](11)--(21);

\draw[->, line width=0.3mm, >=latex, shorten <= 0.2cm, shorten >= 0.15cm](21)--(31);

\draw[->, line width=0.3mm, >=latex, shorten <= 0.2cm, shorten >= 0.15cm](31)--(11);
\draw[->, line width=0.3mm, >=latex, shorten <= 0.2cm, shorten >= 0.15cm](31)--(12);
\end{tikzpicture}
{\caption{$D$.}}
\end{center}
\end{subfigure}%
\begin{subfigure}{0.5\textwidth}
\begin{center}
\tikzstyle{every node}=[circle, draw, fill=black!100,
                       inner sep=0pt, minimum width=5pt]
\begin{tikzpicture}[thick,scale=0.5]%
\draw(0,0)node[label={[yshift=0.1cm]90:{$1$}}](12){};
\draw(0,2)node[label={[yshift=0.1cm]90:{$1$}}](11){};
\draw(-3,-1)node[label={[yshift=0.1cm]90:{$2$}}](21){};
\draw(3,-2)node[label={[yshift=0.1cm]90:{$1$}}](31){};

\draw[->, line width=0.3mm, >=latex, shorten <= 0.2cm, shorten >= 0.15cm](21)--(11);
\draw[dashed, ->, line width=0.3mm, >=latex, shorten <= 0.2cm, shorten >= 0.15cm](12)--(21);

\draw[dashed, ->, line width=0.3mm, >=latex, shorten <= 0.2cm, shorten >= 0.15cm](21)--(31);

\draw[dashed, ->, line width=0.3mm, >=latex, shorten <= 0.2cm, shorten >= 0.15cm](31)--(12);
\draw[->, line width=0.3mm, >=latex, shorten <= 0.2cm, shorten >= 0.15cm](11)--(31);

\end{tikzpicture}
{\caption{$D'$.}}
\end{center}
\end{subfigure}
{\caption{Two tripartite tournaments with score-list parity but are not $TT_3$-equivalent.}\label{figA9.4.7}}
\end{figure}
\end{eg}
\noindent\par Noting that $\mathscr{C}_3=\{C_3, TT_3\}$, Examples \ref{egA9.3.6} and \ref{egA9.4.12} suggest that $\mathscr{C}_3$-reversal is required for the analogue of tripartite tournaments.

\begin{thm}\label{thmA9.4.13}
~\\(i) In the family of tripartite tournaments, $\mathscr{C}_3$ refines $\bigcup\limits_{k\ge 3}\mathscr{C}_k$.
\\(ii) Two tripartite tournaments have score-list parity if and only if they are $\mathscr{C}_3$-equivalent.
\end{thm}
\noindent\textit{Proof}: (i) Let $T$ be a tripartite tournament containing an orientation $F$ of a $k$-cycle, say $u_1 u_2 \ldots u_k u_1$, $k\ge 3$. By Corollary \ref{corA9.4.9}(i), it suffices to consider $k=3$ or $4$, and show that an $F$-reversal in $T$ is equivalent to some sequence of $\mathscr{C}_3$-reversals. The case of $k=3$ is trivial. So, consider $k= 4$. We may assume WLOG that $\{u_1, u_3\} \subseteq V_1$ and $\{u_2, u_4\} \subseteq V_2$. Let $w\in V_3$ and denote the arc joining $u_i$ and $w$ ($u_{i+1}$ resp.) in $T$ as ${a_i}$ ($b_i$ resp.) for $i=1,2,3,4$, with addition modulo $k$. Then, an $F$-reversal is equivalent to the sequence 
\begin{align*}
\{a_1, a_2, b_1\}, \{\tilde{a}_2, a_3, b_2\}, \{\tilde{a}_3, a_4, b_3\}, \{\tilde{a}_4, \tilde{a}_1, b_4\}
\end{align*}
of $\mathscr{C}_3$-reversals.
\indent\par (ii) $(\Rightarrow)$ This follows from Theorem \ref{thmA9.4.4}, (i) and the Refinement Lemma.
\qed

\section{Alternative proofs of Theorem \ref{thmA9.3.1} and Corollary \ref{corA9.3.2}}
In this section, we outline alternative proofs of Theorem \ref{thmA9.3.1} and Corollary \ref{corA9.3.2} based on Theorem \ref{thmA9.1.2} and some elementary ideas in contrast to the involved tools in the previous sections. In addition, we employ a bottom-up approach in which Corollary \ref{corA9.3.2} is proved before Theorem \ref{thmA9.3.1}.
\\
\\\textit{\nth{2} proof of Corollary \ref{corA9.3.2}}: $(\Rightarrow)$ Let $D$ and $D'$ be two orientations of $G(p_1,p_2, \ldots, p_n)$ having the same score list. Since $G$ is bipartite, there exists a bipartition $(X,Y)$ of $G(p_1,p_2,\ldots, p_n)$. If $G(p_1, p_2, \ldots, p_n)\not\cong K(X,Y)$, then add the set of arcs $\overrightarrow{R}=\{xy\mid x\in X, y\in Y, xy\not\in E(G(p_1,p_2, \ldots, p_n))\}$ to $D$ and $D'$ to obtain two orientations $F$ and $F'$, respectively, of $K(X,Y)$ that have the same score list. For convenience, we also denote the underlying edges of $\overrightarrow{R}$ as $R=\{\mathrm{U}(xy)\mid xy\in\overrightarrow{R}\}$. By Theorem \ref{thmA9.1.2}, $F'$ can be obtained from $F$ through a sequence $\mathcal{S}^*=C^{(1)}_4 C^{(2)}_4 \ldots C^{(t)}_4$ of $C_4$-reversals.
\indent\par First, we show that $\mathcal{S}^*$ is equivalent to a sequence $\mathcal{S}$ of $C_{2k}$-reversals, $k\ge 2$, whose underlying cycles are pairwise edge-disjoint; we proceed by induction on the number $t$ of $4$-dicycles in $\mathcal{S}^*$. The base case of $t=1$ is trivial. Suppose $t\ge 2$. By induction hypothesis, the subsequence $\mathcal{S}^*_{sub}=C^{(1)}_4 C^{(2)}_4 \ldots C^{(t-1)}_4$ equivalent to a sequence $\mathcal{S}_{sub}=Z^{(1)} Z^{(2)} \ldots Z^{(r)}$, $1\le r\le t-1$, of $C_{2k}$-reversals where $E(\mathrm{U}(Z^{(i)}))\cap E(\mathrm{U}(Z^{(j)}))=\emptyset$ for any $i\neq j$. Now if $E(\mathrm{U}(C^{(t)}_4))\cap E(\mathrm{U}(Z^{(i)}))=\emptyset$ for all $1\le i\le r$, then we are done. Otherwise, let $Q=E(\mathrm{U}(C^{(t)}_4))\Delta \bigcup\limits_{i=1}^{r} E(\mathrm{U}(Z^{(i)}))$, where $A\Delta B=A\cup B-A\cap B$ represents the symmetric difference of two sets $A$ and $B$. We claim that the arc-induced subdigraph $F^*=\langle Q\rangle_{F}$ consists of only even dicycles, whose underlying cycles are pairwise edge-disjoint. The edge-disjoint property follows from the edge-disjointness of the $Z^{(i)}$'s and the definition of $Q$.
\indent\par
\begin{figure}[h]
\begin{subfigure}{.5\textwidth}
\begin{center}
\tikzstyle{every node}=[circle, draw, fill=black!100,
                       inner sep=0pt, minimum width=5pt]
\begin{tikzpicture}[thick,scale=0.6]%
\draw(0,0)node[label={[yshift=-0.1cm]270:{$a$}}](a){};
\draw(2,0)node[label={[yshift=-0.1cm]270:{$v$}}](v){};
\draw(4,0)node[label={[yshift=-0.1cm]270:{$u$}}](u){};
\draw(6,0)node[label={[yshift=-0.1cm]270:{$b$}}](b){};
\draw(2,2)node[label={[yshift=0.1cm]90:{$x_1$}}](x_1){};
\draw(4,2)node[label={[yshift=0.1cm]90:{$x_2$}}](x_2){};

\draw[dashed, ->, line width=0.3mm, >=latex, shorten <= 0.2cm, shorten >= 0.15cm](-2,0)--(a);
\draw[dashed, ->, line width=0.3mm, >=latex, shorten <= 0.2cm, shorten >= 0.15cm](a)--(v);
\draw[dashed, ->, line width=0.3mm, >=latex, shorten <= 0.2cm, shorten >= 0.15cm](v) to [out=-30,in=-150] (u);
\draw[dashed, ->, line width=0.3mm, >=latex, shorten <= 0.2cm, shorten >= 0.15cm](u)--(b);
\draw[dashed, ->, line width=0.3mm, >=latex, shorten <= 0.2cm, shorten >= 0.15cm](b)--(8,0);

\draw[dotted, ->, line width=0.3mm, >=latex, shorten <= 0.2cm, shorten >= 0.15cm](u) to [out=150,in=30](v);
\draw[dotted, ->, line width=0.3mm, >=latex, shorten <= 0.2cm, shorten >= 0.15cm](v)--(x_1);
\draw[dotted, ->, line width=0.3mm, >=latex, shorten <= 0.2cm, shorten >= 0.15cm](x_1)--(x_2);
\draw[dotted, ->, line width=0.3mm, >=latex, shorten <= 0.2cm, shorten >= 0.15cm](x_2)--(u);
\end{tikzpicture}
{\caption{$Z^{(i)}=\ldots avu b\ldots$ and $C^{(t)}_4=uvx_1 x_2u$.}}
\end{center}
\end{subfigure}%
\begin{subfigure}{0.5\textwidth}
\begin{center}
\tikzstyle{every node}=[circle, draw, fill=black!100,
                       inner sep=0pt, minimum width=5pt]
\begin{tikzpicture}[thick,scale=0.6]%
\draw(0,0)node[label={[yshift=-0.1cm]270:{$a$}}](a){};
\draw(2,0)node[label={[yshift=-0.1cm]270:{$v$}}](v){};
\draw(4,0)node[label={[yshift=-0.1cm]270:{$u$}}](u){};
\draw(6,0)node[label={[yshift=-0.1cm]270:{$b$}}](b){};
\draw(2,2)node[label={[yshift=0.1cm]90:{$x_1$}}](x_1){};
\draw(4,2)node[label={[yshift=0.1cm]90:{$x_2$}}](x_2){};

\draw[->, line width=0.3mm, >=latex, shorten <= 0.2cm, shorten >= 0.15cm](-2,0)--(a);
\draw[->, line width=0.3mm, >=latex, shorten <= 0.2cm, shorten >= 0.15cm](a)--(v);
\draw[->, line width=0.3mm, >=latex, shorten <= 0.2cm, shorten >= 0.15cm](u)--(b);
\draw[->, line width=0.3mm, >=latex, shorten <= 0.2cm, shorten >= 0.15cm](b)--(8,0);

\draw[->, line width=0.3mm, >=latex, shorten <= 0.2cm, shorten >= 0.15cm](v)--(x_1);
\draw[->, line width=0.3mm, >=latex, shorten <= 0.2cm, shorten >= 0.15cm](x_1)--(x_2);
\draw[->, line width=0.3mm, >=latex, shorten <= 0.2cm, shorten >= 0.15cm](x_2)--(u);
\end{tikzpicture}
{\caption{$Z^{(i)}_{new}=\ldots avx_1 x_2 u b\ldots$.}}
\end{center}
\end{subfigure}
{\caption{Combining $Z^{(i)}$ and $C^{(t)}_4$ to become $Z^{(i)}_{new}$.}\label{figA9.5.8}}
\end{figure}
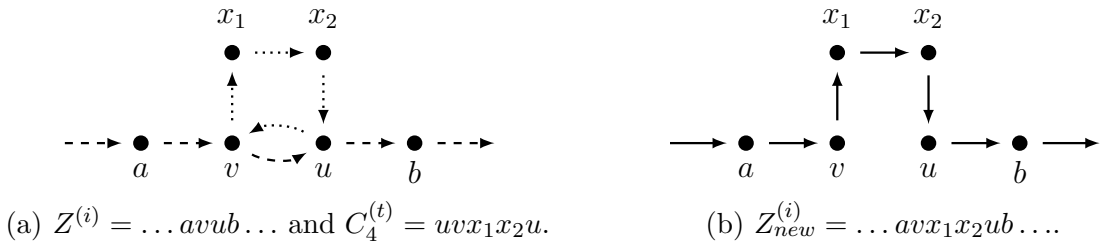

 To deduce that $F$ contains only dicycles, note if $uv\in E(\mathrm{U}(C^{(t)}_4))\cap E(\mathrm{U}(Z^{(i)}))$ for some $1\le i\le r$, we must have $uv$ oriented differently in $Z^{(i)}$ and $C^{(t)}_4$, say $Z^{(i)}=\ldots avu b\ldots$ and $C^{(t)}_4=uvx_1 x_2u$, so that the latter is a feasible dicycle reversal; and they combine to give (and will be replaced by) the dicycle $Z^{(i)}_{new}=\ldots avx_1 x_2 u b\ldots$ (see Figure \ref{figA9.5.8}). If there exists some edge $xy\in E(\mathrm{U}(Z^{(i)}_{new}))\cap E(\mathrm{U}(C^{(j)}_4))$ (certainly, $xy\neq uv$) for some $1\le j\le r$, $j\neq i$, then repeat the argument on $Z^{(i)}_{new}$ and $C^{(j)}_4$ in place of $C^{(t)}_4$ and $C^{(i)}_4$ respectively. If there exists some arcs $xy, yx\in A(Z^{(i)}_{new})$, then we deduce via similar considerations that $Z^{(i)}_{new}$ can be replaced by at most two new dicycles. So,  $\mathcal{S}^*$ is equivalent to a sequence $\mathcal{S}$ of $C_{2k}$-reversals; these dicycles in $F^*$ must be of even length since the underlying graph of $F^*$ is a subgraph of $K(X,Y)$. This completes the induction.
\indent\par Lastly, we show that none of the underlying cycles $\mathrm{U}(Z^{(i)})$'s contains an edge in $R$. Since any arc in $\overrightarrow{R}$ is oriented the same in $F$ and $F'$, the number of times an arc $uv$ appears in the dicycles $Z^{(i)}$'s is the same as that of its converse $vu$, that is $|\{i \mid uv\in A(Z^{(i)})\}|=|\{i \mid vu\in A(Z^{(i)})\}|$. By the edge-disjoint property of $Z^{(i)}$'s in $\mathcal{S}$, it follows that $A(Z^{(i)})\cap\overrightarrow{R}=\emptyset$ for all $i$. Hence, $\mathcal{S}$ is a feasible sequence of $C_{2k}$-reversals to obtain $D'$ from $D$.
\qed

\indent\par For our next proof, we require the subdivision operation. The \textit{subdivision} of an edge (arc resp.) $xy$ in a graph (digraph resp.) is the operation of removing $xy$ and adding a new vertex $w_{xy}$ with edges (arcs resp.) $xw_{xy}$ and $w_{xy}y$.
\\\\
\noindent\textit{\nth{2} proof of Theorem \ref{thmA9.3.1}}: $(\Rightarrow)$ Let $D$ and $D'$ be two orientations of $G(p_1,p_2, \ldots, p_n)$ having the same score list and $H$ be the bipartite graph obtained by subdividing every edge of $G$. Suppose $1, 2, \ldots, n$ is a fixed labelling of the vertices in $G$. Let $F$ and $F'$ be the orientations of $H(p_1, p_2, \ldots, p_n, q_1, q_2, \ldots, q_m)$ obtained from the subdivision of every arc of $D$ and $D'$ respectively; and the addition of the arcs $uw_{xy}$ and $w_{xy}v$ for any $u\in V_i-\{x,y\}$ and $v\in V_j-\{x,y\}$ with $i<j$, if $xy$ is an arc between the partite sets $V_i$ and $V_j$. Formally,
\begin{align*}
V(F)=\ &V(D)\cup \{w_{xy} \mid xy\in A(D)\}\text{ and}\\
A(F)=\ &\{xw_{xy}, w_{xy}y \mid xy\in A(D)\} \cup \overrightarrow{R},\text{ where}\\
\overrightarrow{R}=\ &\{uw_{xy}, w_{xy}v\mid u\in V_i-\{x,y\}, v\in V_j-\{x,y\}, i<j, xy\in A(D)\\
\ &\text{and } x\in V_i, y\in V_j\text{ or } x\in V_j, y\in V_i\}.
\end{align*}
The same holds for $F'$ and $D'$.
\noindent\par We verify that $F$ and $F'$ have the same score list as follows. Let $ij\in E(G)$ with $i<j$. Note that the set of vertices added from the subdivisions of the arcs between partite sets $V_i$ and $V_j$, namely $Q_{ij}=\{w_{xy}\mid x\in V_i, y\in V_j \text{ or }x\in V_j, y\in V_i\}$, forms a partite set in $H(p_1, p_2, \ldots, p_n, q_1, q_2, \ldots, q_m)$ with size $|Q_{ij}|=q_{ij}=|V_i||V_j|$. Clearly, each vertex $w\in Q_{ij}$ has score $|V_j|$ in $F$ and $F'$. Furthermore, the added arcs incident to $Q_{ij}$ increases the score of each $u\in V_i$ by $(|V_i|-1)|V_j|$ and do not affect the score of each $v\in V_j$.
\indent\par By Corollary \ref{corA9.3.2}, there exists some sequence $\mathcal{S}$ of $C_{2k}$-reversals, $k\ge 2$, whose underlying graphs are pairwise edge-disjoint, to obtain $F'$ from $F$. Since $\overrightarrow{R}\subset A(F)\cap A(F')$, it follows from the edge-disjoint property that none of the arcs in $\overrightarrow{R}$ are used in $\mathcal{S}$. Hence, restoring any subdivided arcs used in $\mathcal{S}$ gives a feasible sequence of $C_{k}$-reversals $k\ge 3$, to obtain $D'$ from $D$.
\qed
\section{Conclusion}
Let us conclude the paper with some possible directions for future research. In this paper, we focus on determining the $C_k$-classes and $\mathscr{C}_k$-classes, $k\ge 3$, in families of orientations of a $G$ vertex-multiplication to extend Theorems \ref{thmA9.1.1} and \ref{thmA9.1.2}. Other digraphs such as the dipaths, antidirected dipaths and dicycles are also studied on tournaments by Reid \cite{KBR}, Rosenfeld \cite{MR1, MR2} and Waldrop \cite{CW1}. We state some of these well-known results. Let $P_{k+1}$ denote a $k$-dipath and $\mathscr{P}_{k+1}$ denote the family of all orientations of a path of length $k$, i.e., $|V(P_{k+1})|=|V(\mathscr{P}_{k+1})|=k+1$. An \textit{antidirected} $k$-dipath (antidirected $k$-dicycle resp.) has every two adjacent arcs in opposite directions (except possibly the first and last arc resp.). Obviously if $k$ is even, then there are two antidirected $k$-dipaths which are converses of each other; we denote the one that starts with a forward arc as $AP_{k+1}$. If $k$ is odd, then the antidirected $k$-dipath is unique (up to isomorphism). Let $AC_k$ denote an antidirected $k$-dicycle.

\begin{thm} (Reid \cite{KBR})
Any two $n$-tournaments are $P_k$-equivalent for each integer $2\le k\le n$.
\end{thm}
\begin{thm} (Rosenfeld \cite{MR1})
Every vertex in a tournament $T$ of order $n\ge 9$ is an end-vertex of some antidirected Hamiltionian path in $T$.
\end{thm}
\begin{thm} (Rosenfeld \cite{MR2})
If $T$ is a tournament of even order $n\ge 28$, then $T$ has an antidirected Hamiltonian cycle.
\end{thm}
\begin{thm} (Waldrop \cite{CW1})
If $n\ge 11$, then any two tournaments of order $n$ are $AP_k$-equivalent for each integer $2\le k\le n-1$.
\end{thm}

\noindent\par Therefore, it is natural to seek analogues for orientations of $G$ vertex-multiplications of the above results. Following Waldrop \cite{CW1} and in view of the structure of $G$ vertex-multiplications, the digraph $M_k$ consisting of $k$ vertex-disjoint arcs from one partite set to another seems interesting for investigation too. Formally, $V(M_k)=X\cup Y$, where $X=\{x_1, x_2, \ldots, x_k\}\subseteq V_x$, $Y=\{y_1, y_2, \ldots y_k\}\subseteq V_y$ for some $xy\in E(G)$ and $A(M_k)=\{x_i y_i\mid i=1,2, \ldots, k\}$. We end by proposing the following problem.

\begin{prob}
Let $G$ be a graph. In the family of orientations of $G(p_1, p_2, \ldots, p_n)$, determine the $\mathscr{F}$-classes where
\\(1) $\mathscr{F}=\{P_{k+1}\}$ (i.e., $k$-dipaths);
\\(2) $\mathscr{F}=\mathscr{P}_{k+1}$ (i.e., family of orientations of $k$-path);
\\(3) $\mathscr{F}=\{AP_{2k+1}, \tilde{AP}_{2k+1}, AP_{2k}\}$ (i.e., antidirected dipaths);
\\(4) $\mathscr{F}=\{AC_k\}$ (i.e., antidirected dicycles);
\\(5) $\mathscr{F}=\{M_k\}$ (i.e., $k$ vertex-disjoint arcs from one partite set to another).
\end{prob}

\section*{Acknowledgements}
The first author would like to thank the National Institute of Education, Nanyang Technological University of Singapore, for the generous support of the Nanyang Technological University Research Scholarship.

\end{document}